\documentclass[11pt]{article}

\usepackage[latin1]{inputenc}

\usepackage{amsfonts,amsmath,amssymb,amsthm,graphicx,epsfig,float}
\usepackage[T1]{fontenc}

\usepackage[english,francais]{babel}

{
  \newtheorem{theoreme}{Th\'eor\`eme}
  \newtheorem*{theoreme*}{Th\'eor\`eme}
  \newtheorem{lemme}[theoreme]{Lemme}

  \newtheorem{definition}{D\'efinition}
  \newtheorem{proposition}[theoreme]{Proposition}
\newtheorem*{corollaire*}{Corollaire}
\newtheorem*{proposition*}{Proposition}
\theoremstyle{remark}
  \newtheorem*{remarque*}{Remarque}
}

\newcounter{ex}

\newenvironment{rem*}{
  \noindent\textbf{Remarque. }}{}

\newcommand{\Cc}{\mathbb{C}}
\newcommand{\Nn}{\mathbb{N}}

\newcommand{\Zz}{\mathbb{Z}}
\newcommand{\Pp}{\mathbb{P}}

\newcommand{\Bcal}{\mathcal{B}}

\newcommand{\Dcal}{\mathcal{D}}
\newcommand{\Ecal}{\mathcal{E}}
\newcommand{\Gcal}{\mathcal{G}}

\newcommand{\Qcal}{\mathcal{Q}}
\newcommand{\Scal}{\mathcal{S}}

\title{{\bf Sur la construction de mesures selles}}
\author{Henry de Thélin}
\date{}

\begin{document}
\maketitle


\def\figurename{{Fig.}}%
\def\proofname{Preuve}
\def\contentsname{Sommaire}%

\begin{abstract}

Nous construisons des mesures selles (dans un sens faible) pour les
endomorphismes holomorphes de $\Pp^2(\Cc)$.

\end{abstract}

\selectlanguage{english}
\begin{center}
{\bf{On saddle measures }}
\end{center}

\begin{abstract}

We build saddle measures (in a weak sense) for holomorphic
endomorphisms of $\Cc \Pp^2$.

\end{abstract}

\selectlanguage{francais}

Mots-clefs: dynamique holomorphe, entropie, exposants de Lyapunov.\\
Classification: 32H50, 37FXX.

\section*{{\bf Introduction}}
\par

A partir d'un endomorphisme holomorphe de $\Pp^2(\Cc)$, $f$, de degré
$d \geq 2$, J.E. Forn{\ae}ss et N. Sibony ont défini le courant de
Green $T$ associé à $f$ (voir \cite{FS} et \cite{FS1}), dont le
support est l'ensemble de Julia de $f$. Ce courant possède un
potentiel continu: on peut donc définir son auto-intersection $\mu= T
\wedge T$ (voir \cite{FS}). La mesure $\mu$ ainsi obtenue est l'unique
mesure d'entropie maximale $2 \log(d)$ (voir \cite{BD2}) et elle a ses
exposants de Lyapunov minorés par $\frac{\log(d)}{2}$ (voir
\cite{BD1}).

L'objet de cet article est de décrire la dynamique de $f$ en dehors du
support de $\mu$. L'entropie topologique de $f$ hors de ce support
étant majorée par $\log(d)$ (voir \cite{Det2}), il s'agira d'une part
de construire des mesures d'entropie $\log(d)$ et d'autre part, d'évaluer leurs
exposants de Lyapunov.

Quand $f$ est un endomorphisme holomorphe hyperbolique (dans un sens
fort), cela a été réalisé par J.E. Forn{\ae}ss et
N. Sibony dans $\cite{FS3}$. En effet, soient $L$ une droite
projective de $\Pp^2(\Cc)$ et $S$ une valeur d'adhérence de $S_m=\frac{1}{m} \sum_{i=0}^{m-1}
\frac{[f^i(L)]}{d^i}$. Le courant $S$ vérifie $f_{*}S=d S$ et en le tranchant avec $T$, on obtient une
mesure $\nu=T \wedge S$ invariante par $f$. Alors, dans \cite{FS3},
J.E. Forn{\ae}ss et N. Sibony ont démontré entre autres que ces mesures
étaient selles (i.e. qu'elles ont un exposant de Lyapunov strictement
positif et un strictement négatif).

Quand $f$ est un endomorphisme holomorphe quelconque de $\Pp^2(\Cc)$,
nous pouvons construire comme précédemment des mesures $\nu=T \wedge
S$ invariantes par $f$. L'objectif de cet article est de voir que ces
mesures permettent de décrire la dynamique de $f$ en dehors du support
de $\mu$. Plus précisément, nous aurons tout d'abord (comme dans \cite{BS3}) le

\begin{theoreme}{\label{entropie}}
L'entropie métrique de $\nu$ est minorée par $\log(d)$.

\end{theoreme}

En particulier, quand le support de la mesure $\nu$ est disjoint de celui
de la mesure $\mu$, la mesure $\nu$ est d'entropie maximale dans le
complémentaire du support de $\mu$ (voir \cite{Det2}).

Les mesures $\nu$ ne sont pas ergodiques en général (par exemple si
$f([z:w:t])=[z^2:w^2:t^2]$, elles ont en général trois
composantes ergodiques). En particulier, le théorème précédent associé
à l'inégalité de Ruelle (voir \cite{Ru} ou le paragraphe \ref{exposant}) ne permet pas d'obtenir que le plus grand
exposant de Lyapunov de $\nu$ est supérieur à $0$ en presque tout
point. Cependant, nous atteindrons ce résultat en
adaptant les arguments de R. Dujardin (voir \cite{Du3}) à notre
situation. Nous aurons donc le

\begin{theoreme}{\label{th2}}

Pour $\nu$ presque tout point $x$, le plus grand exposant de Lyapunov
en $x$ est supérieur ou égal à $\frac{\log(d)}{2}$.

\end{theoreme}

Pour décrire la dynamique de $f$ hors du support de $\mu$, il nous reste à
estimer le plus petit exposant de Lyapunov de $\nu$ hors de ce
support. C'est l'objet du

\begin{theoreme}{\label{th3}}
Supposons que $\nu$ ne charge pas les courbes algébriques.

Alors, pour $\nu$ presque tout point $x$ hors du support de $\mu$, le plus
petit exposant de Lyapunov $\chi_1(x)$ est négatif ou nul.

\end{theoreme}

Les mesures $\nu$ sont donc faiblement selles en dehors du support de
$\mu$.

Par ailleurs, dans le paragraphe \ref{exemple},
nous donnerons un exemple de mesure $\nu$ qui ne charge aucune courbe
algébrique pour laquelle le plus petit exposant de Lyapunov est nul.

Remarquons que si $\nu=T \wedge S$ charge une courbe algébrique alors
celle-ci est nécessairement prépériodique. En effet, d'une part $\nu$
est invariante et d'autre part, comme $T$ est à potentiel höldérien, $\nu$ ne charge aucun point. En particulier, nous avons le

\begin{corollaire*}

Si $f$ ne possède aucune courbe périodique alors toutes les mesures
$\nu=T \wedge S$ construites précédemment ont leur plus petit exposant
de Lyapunov négatif ou nul pour presque tout point hors du support de $\mu$.

\end{corollaire*}

Voici maintenant le plan de ce texte: dans le premier paragraphe, nous minorerons
l'entropie de la mesure $\nu$ par $\log(d)$. La seconde et la
troisième partie de cet article seront consacrées à des rappels d'une
part sur la théorie de Pesin et d'autre part sur les courants tissés
(ou géométriques) introduits dans \cite{Di}. Ensuite, le quatrième
paragraphe démontrera la minoration du plus grand exposant de Lyapunov
de $\nu$, tandis que le cinquième traitera de la majoration du plus
petit. Enfin, dans la sixième partie, nous donnerons l'exemple de
mesure $\nu$ qui rend le théorème \ref{th3} optimal.

$ $

{\bf Remerciements:} Je tiens à remercier T.-C. Dinh et R. Dujardin
pour leurs remarques sur le fond et la forme de cet article.

\section{{\bf Entropie métrique de la mesure $\nu$}}

Dans ce paragraphe, on considère $L$ une droite projective et $S$ une
valeur d'adhérence de $\frac{1}{n} \sum_{i=0}^{n-1}
\frac{[f^i(L)]}{d^i}$. Le courant $S$ vérifie $f_{*}S=dS$. Par
ailleurs, en tranchant ce courant avec $T$, on obtient une mesure
$\nu=T \wedge S$ qui est invariante par $f$.

L'objectif de ce paragraphe est alors de démontrer le

\begin{theoreme*}
L'entropie métrique de $\nu$ est minorée par $\log(d)$.

\end{theoreme*}

Avant de passer à la démonstration de ce théorème, nous allons faire
quelques rappels sur l'entropie.

\subsection{{\bf Entropie métrique}}

Notons $d_n(x,y) = \max_{0 \leq i \leq n-1} \{ d(f^i(x), f^i(y)) \}$
et $B_n(x, \epsilon)$ la boule de centre $x$ et de rayon $\epsilon$
pour cette métrique.

A partir du théorème de Brin-Katok (voir \cite{BK}) nous pouvons
définir l'entropie métrique $h_{\nu}(f)$ de $\nu$ par:
$$h_{\nu}(f)= \int  \sup_{\epsilon > 0} \liminf_{n} - \frac{1}{n} \log
(\nu(B_n(x, \epsilon))) d\nu(x).$$
Cette quantité décrit donc la décroissance moyenne de la masse d'une
boule dynamique $B_n(x, \epsilon)$ pour $\nu$.

Dans \cite{BS3}, E. Bedford et J. Smillie ont donné une méthode pour
minorer l'entropie métrique de certaines mesures. Elle s'appuie sur la
démonstration du principe variationnel (voir \cite{W}) et plus
précisément sur la

\begin{proposition*}

Fixons $\epsilon >0$. Soient $\sigma_n$ une suite de probabilités et
$\nu_n = \frac{1}{n} \sum_{i=0}^{n-1} f_{*}^i \sigma_n$. Si $\nu_{n_j}$ converge vers $\nu$ et
$\sigma_{n_j}(B_{n_j}(x,\epsilon)) \leq c_{n_j}$ pour toute boule
dynamique $B_{n_j}(x,\epsilon)$ alors:
$$h_{\nu}(f) \geq \limsup - \frac{1}{n_j} \log c_{n_j}.$$

\end{proposition*}

\subsection{{\bf Minoration de l'entropie de $\nu$}}

La démonstration du théorème \ref{entropie} va se faire en deux étapes. Dans la
première, on va utiliser la proposition précédente avec $\sigma_n=
\frac{f^{n*} \omega}{d^n} \wedge [L]$ (où $\omega$ est la forme de
Fubini-Study de $\Pp^2(\Cc)$). Celle-ci nous donnera alors,
exactement comme dans \cite{BS3}, une mesure $\nu^{'}$ d'entropie
minorée par $\log(d)$. Il restera alors à voir que $\nu$ et $\nu^{'}$
sont égales.

\subsubsection{Construction d'une mesure d'entropie minorée par
  $\log(d)$}

Quitte à extraire une sous-suite, nous supposerons que
$S_n=\frac{1}{n} \sum_{i=0}^{n-1} \frac{[f^i(L)]}{d^i}$ converge
vers $S$. La preuve étant la même que celle de E. Bedford et
J. Smillie, nous passerons vite sur certains points.

Appliquons la proposition précédente à  $\sigma_n= \frac{f^{n*} \omega}{d^n} \wedge [L]$. Nous avons
$$ \sigma_{n} ( B_{n}(x,\epsilon)) = \frac{1}{d^{n}}
  \int_{B_{n}(x,\epsilon)} f^{n*} \omega \wedge [L] \leq
  \frac{1}{d^{n}} v^0(f,n,\epsilon)$$
où $v^0(f,n,\epsilon)$ est le suprémum sur toutes les boules
  dynamiques $B_{n}(x,\epsilon)$ du volume de $f^{n}(B_{n}(x,\epsilon)
  \cap L)$ compté avec multiplicité. Si $\nu_n^{'}$ désigne la mesure $\frac{1}{n} \sum_{i=0}^{n-1} f_{*}^i
  \sigma_n$ et $\nu_{n_j}^{'}$ une sous-suite de 
$\nu_n^{'}$ qui converge vers $\nu{'}$, on a:
$$h_{\nu^{'}}(f) \geq \limsup (- \frac{1}{n_j} \log (
  \frac{1}{d^{n_j}}v^0(f,n_j,\epsilon))).$$
En utilisant le théorème de Yomdin (voir \cite{Y}) on obtient alors la
  minoration cherchée.

Il reste à voir que $\nu$ et $\nu^{'}$ sont égales. C'est l'objet du
  paragraphe suivant.

\subsubsection{Minoration de l'entropie de $\nu$}

Commençons par rappeler brièvement la construction de $T$ (voir
\cite{FS} et \cite{FS1}).

Comme la forme $f^{*}\omega$ est cohomologue à $d \omega$, on a:
$$\frac{f^{*}\omega}{d}= \omega +dd^c u,$$
où $u$ est une fonction lisse de $\Pp^2(\Cc)$.

En itérant cette relation, on obtient que:
$$\frac{f^{i*}\omega}{d^i}= \omega +dd^c G_i,$$
avec $G_i=\sum_{l=0}^{i-1} \frac{u \circ f^l}{d^l}$.

Autrement dit, en passant à la limite, on a $T= \omega + dd^c G$ où $G$ est une fontion continue qui vérifie:
$$\max_{\Pp^2(\Cc)}|G_{i}-G| \leq \frac{C}{d^i}.$$
Maintenant
$$\nu_{n_j}^{'}=\frac{1}{n_j} \sum_{i=0}^{n_j -1} f_{*}^i( \frac{f^{n_j*}
  \omega}{d^{n_j}} \wedge [L])=\frac{1}{n_j} \sum_{i=0}^{n_j -1}
  \frac{f^{n_j -i*}
  \omega}{d^{n_j -i}} \wedge \frac{f_{*}^i([L])}{d^i}$$
qui est égal à 
$$\frac{1}{n_j} \sum_{i=0}^{n_j -1}  (\frac{f^{n_j -i*}
  \omega}{d^{n_j -i}}-T) \wedge \frac{f_{*}^i([L])}{d^i} + \frac{1}{n_j} \sum_{i=0}^{n_j -1}
  T \wedge \frac{f_{*}^i([L])}{d^i}.$$
Le dernier terme est $T \wedge S_{n_j}$ et converge donc vers $\nu$
  car $T$ est à potentiel continu.

Par ailleurs, en appliquant au premier terme une fonction test $\Phi$
  et en utilisant la construction de $T$ on obtient:
$$\int \frac{1}{n_j} \sum_{i=0}^{n_j -1}  (G_{n_j -i} - G) dd^c \Phi
  \wedge \frac{f_{*}^i([L])}{d^i}$$
qui en valeur absolue est majoré par $|\Phi|_{C^2} \frac{1}{n_j}
  \sum_{i=0}^{n_j -1} \frac{C}{d^{n_j - i}}$ qui tend bien vers $0$
  quand $j$ croît vers l'infini.

Cela montre bien que $\nu$ et $\nu^{'}$ sont égales.

\section{{\bf Un peu de théorie de Pesin}}{\label{exposant}}

Dans tout cet article, $\Pp^2(\Cc)$ sera muni d'une famille de cartes
holomorphes $\tau_x : \Cc^2 \mapsto \Pp^2(\Cc)$ avec $\tau_x(0)=x$
et $(x,z) \mapsto \tau_x(z)$ localement $C^{\infty}$.

Soient $f$ un endomorphisme holomorphe de $\Pp^2(\Cc)$ et $\nu$ une
mesure de probabilité invariante.

Pour faire de la théorie de Pesin associée à $f$, nous allons définir
l'extension naturelle $\widehat{\Pp^2(\Cc)}$ de $(\Pp^2(\Cc),f,\nu)$: c'est l'ensemble
$\widehat{\Pp^2(\Cc)}= \{ \widehat{x}=(...,x_{-n},...,x_0) \in \Pp^2(\Cc)^{\Zz^{-}}
\mbox{  ,  } f(x_{-n})=x_{-n+1} \}$, des préhistoires des points de
$\Pp^2(\Cc)$. Dans cet espace $f$ induit une application $\widehat{f}$
qui est le décalage à droite et si on note $\pi$ la projection
canonique $\pi(\widehat{x})=x_{0}$, alors $\nu$ se relève par $\pi$ en
une unique mesure de probabilité $\widehat{\nu}$ invariante par $\widehat{f}$ qui
vérifie $\pi_{*} \widehat{\nu}= \nu$.

Maintenant, à partir de $f_x = \tau_{f(x)}^{-1} \circ f \circ \tau_x$,
on peut définir une application $D \widehat{f}(\widehat{x})$ de
$\widehat{X}$ à valeur dans $\mbox{Mat}(2, \Cc)$ en posant: $D
\widehat{f}(\widehat{x})= Df_x (0)$ (où $\pi(\widehat{x})=x$).

C'est à  $D \widehat{f}$ que l'on va appliquer une extension du
théorème d'Osedelec (qui est valable dans un cadre non intégrable) et
la théorie de Pesin.

\begin{theoreme*}{\label{exposants}}(Voir \cite{PS} p.35).
Il existe un borélien invariant $\widehat{X} \subset \widehat{\Pp^2(\Cc)}$ avec
$\widehat{\nu}(\widehat{\Pp^2(\Cc)} - \widehat{X})=0$ tel que pour tout
  $\widehat{x}$ dans $\widehat{X}$ on ait:

1) Une décomposition mesurable de $\Cc^2$ en espaces complexes de la forme:
$$\Cc^2= \oplus_{i=1}^{k(\widehat{x})} E_i(\widehat{x})$$
qui vérifie $D \widehat{f}(\widehat{x})
E_i(\widehat{x})=E_i(\widehat{f} (\widehat{x}) )$.

2) L'existence de fonctions invariantes par $\widehat{f}$ (qui peuvent
valoir $- \infty$):
$$\chi_1(\widehat{x}) < \cdots <\chi_{k(\widehat{x})}(\widehat{x}),$$
avec
$$ \lim_{m \rightarrow + \infty} \frac{1}{m} \log
\| A(\widehat{x},m)v \|=  \chi_i(\widehat{x})$$
pour tout $v \in E_i(\widehat{x}) \setminus \{ 0\}$.

Ici les $A(\widehat{x},m)$ sont définis par les relations:
$$A(\widehat{x},m)= D \widehat{f}(\widehat{f}^{m-1}(\widehat{x})) \circ
  \cdots \circ  D \widehat{f}(\widehat{x}) \mbox{       pour } m > 0.$$
3) Si $k(\widehat{x})=2$ alors:
$$\lim_{m \rightarrow  \infty} \frac{1}{m} \log \sin ( \angle (E_1
(\widehat{f}^{m}(\widehat{x})),E_2
(\widehat{f}^{m}(\widehat{x}))))=0.$$

\end{theoreme*}

Les réels $\chi_i (\widehat{x})$ sont les exposants de Lyapunov de
$\widehat{f}$. Par ailleurs, si on a
$\pi(\widehat{x})=\pi(\widehat{y})=x$ avec $\widehat{x}$ et
$\widehat{y}$ dans $\widehat{X}$, on a $k(\widehat{x})=k(\widehat{y})$
et $\chi_i (\widehat{x})=\chi_i (\widehat{y})$. Cela permet donc de définir
les exposants de Lyapunov pour $\nu$-presque tout point $x$.

Dans la suite, nous nous plaçons sur $\widehat{X}^{'}$ qui est
l'ensemble des points de $\widehat{X}$ pour lesquels $k(\widehat{x})=2$ (c'est-à-dire $\widehat{\nu}$ a deux exposants de Lyapunov distincts) et
$\chi_1(\widehat{x}) > 0$. C'est un borélien invariant par
$\widehat{f}$. En particulier, si on considère la mesure
$\widehat{\nu}^{'}$ définie par $\widehat{\nu}^{'}(A)=\widehat{\nu}(A
\cap \widehat{X}^{'})$, on obtient une mesure invariante par
$\widehat{f}$ pour laquelle
$\log ^{+} \| (D \widehat{f})^{-1} \|$ est intégrable (voir
l'appendice $A.1$ de \cite{K}).

Nous pouvons donc lui appliquer le théorème de $\gamma$-réduction de Pesin
(voir \cite{KH}) et ainsi obtenir:

\begin{theoreme*}($\gamma$-réduction de Pesin (voir \cite{KH})).

Pour tout $\gamma > 0$, il existe une application $C_{\gamma} :
\widehat{X}^{'} \mapsto GL(2,\Cc)$ et un borélien invariant
$\widehat{Y}^{'} \subset \widehat{X}^{'}$ avec
$\widehat{\nu}(\widehat{X}^{'} - \widehat{Y}^{'})=0$ tel que pour tout
$\widehat{x} \in \widehat{Y}^{'}$ on ait:

1) $$\lim_{m \rightarrow  \infty} \frac{1}{m} \log \| C_{\gamma}^{\pm 1} (
\widehat{f}^{m}(\widehat{x})) \|=0$$
(on parle de fonction tempérée).

2) La matrice $A_{\gamma}(\widehat{x})= C_{\gamma}^{-1} (
\widehat{f}(\widehat{x})) D \widehat{f} (\widehat{x})
C_{\gamma} (\widehat{x})$ est égale à $\left( \begin{array}{cc}
\mu_1(\widehat{x}) & 0\\
0 & \mu_2(\widehat{x})\\
\end{array} \right) $
avec 
$$e^{\chi_i (\widehat{x})- \gamma} \leq |\mu_i(\widehat{x})|
\leq e^{\chi_i (\widehat{x})+ \gamma}.$$
3) Enfin, $C_{\gamma} (\widehat{x})$ envoie la décomposition
standard de $\Cc^2$ sur $E_1(\widehat{x}) \oplus E_2(\widehat{x})$.

\end{theoreme*}

Maintenant, si on note $g_{\widehat{x}}$ la lecture de $f_x$ dans ces cartes
(i.e.  $g_{\widehat{x}}= C_{\gamma}^{-1} (\widehat{f}(\widehat{x})) \circ f_x \circ
C_{\gamma} (\widehat{x})$ où $\pi(\widehat{x})=x$), on en
déduit la proposition suivante (voir \cite{KH} p. 673):

\begin{proposition*}
Pour tout $\widehat{x} \in \widehat{Y}^{'}$ les $g_{\widehat{x}}$ sont définies sur des boules
$B(\delta(\widehat{x}))$ où $\delta(\widehat{x})$ sont des fonctions
tempérées.

Par ailleurs, on a $g_{\widehat{x}}(0)=0$ et $Dg_{\widehat{x}}(0)=\left( \begin{array}{cc}
\mu_1(\widehat{x}) & 0\\
0 & \mu_2(\widehat{x})\\
\end{array} \right)$. Enfin si $g_{\widehat{x}}(w)= Dg_{\widehat{x}}(0)w
+h(w)$ alors $\| D_w h \| \leq \frac{ \| w \| \gamma}{
  \delta(\widehat{x}) }$ pour $w \in B(\delta(\widehat{x}))$.

\end{proposition*}

\begin{remarque*}

Comme $\delta(\widehat{x})$ est tempérée, nous pouvons supposer dans
la proposition précédente que
$e^{- \gamma} <
\frac{\delta(\widehat{f}(\widehat{x}))}{\delta(\widehat{x})} <
e^{\gamma}$ (voir \cite{KH} p. 668).
 \end{remarque*}

Signalons enfin que les exposants de Lyapunov sont reliés à l'entropie
par le

\begin{theoreme*} Inégalité de Ruelle, cas complexe (voir \cite{Ru}).

Soient $f$ un endomorphisme holomorphe de $\Pp^2$ et $\nu$ une
mesure invariante par $f$. On  a alors:
$$\frac{h_{\nu}(f)}{2} \leq \int \chi^{+} (x) d\nu(x)$$
où $\chi^{+} (x)= \sum_{i \mbox{, } \chi_i (x) > 0}  \chi_i (x) .$

\end{theoreme*}

\section{{\bf Courants géométriques}}{\label{courant}}

Dans \cite{Di}, T.-C. Dinh a démontré que le courant $S$ est
tissé. Cela signifie que son support contient beaucoup de disques
analytiques.

Dans ce paragraphe, nous allons tout d'abord rappeler la notion de
courant tissé (que nous appellerons courant géométrique), puis nous étudierons
l'intersection de $S$ avec $T$, en suivant de près l'article de
R. Dujardin (voir \cite{Du2}).

\subsection{{\bf Courants géométriques}}

Considérons un ouvert $\Omega$ de $\Cc^2$ et $T$ un $(1,1)$ courant
positif.

\begin{definition}

Le courant $T$ est uniformément géométrique dans $\Omega$ si pour tout
$x \in \mbox{Supp}(T) \cap \Omega$, il existe un bidisque $B$, un
ouvert $U$ de $B$ contenant $x$ et une constante $c(B)$ tels que:
$$T_{|U}= \int_{\Gamma \in \mathcal{G}} [\Gamma \cap U] d \lambda
(\Gamma).$$
Ici $ \mathcal{G}$ est l'ensemble des sous-ensembles analytiques de
$B$ de masse inférieure à $c(B)$ et $\lambda$ est une mesure sur cet
espace compact.

\end{definition}

Remarquons que les $\Gamma$ ne sont pas supposés disjoints.

De façon analogue aux courants laminaires (voir par exemple \cite{C},
\cite{Det1}, \cite{Di} et \cite{Du1}), on peut définir:

\begin{definition}

Un courant $T$ est géométrique dans $\Omega$ s'il existe une suite
d'ouverts $\Omega_i \subset \Omega$ avec $||T||(\partial \Omega_i)=0$
et une suite croissante $(T_i)_{i \geq 0}$, $\mbox{ }T_i$ uniformément
géométrique dans $\Omega_i$ tels que $\lim_{i \rightarrow \infty}
T_i=T$.

\end{definition}

\subsection{{\bf Caractère géométrique de $S$}}{\label{geometrique}}

Dans notre contexte, $S$ est un courant géométrique (voir \cite{Di}).
Cependant, nous allons voir que de façon analogue à R. Dujardin pour
le cas laminaire, nous pouvons raffiner la suite de courants
uniformément géométriques $S_k$ qui
croît vers $S$ de sorte que $T \wedge S_k$ croisse aussi vers $T
\wedge S$.

On commence par fixer un ouvert $\Omega$ et deux projections linéaires
génériques $\pi_1$ et $\pi_2$. On peut alors recouvrir
$\overline{\Omega}$ par une subdivision en $4$-cubes affines:
$$\Qcal = \{ \pi_1^{-1}(s_1) \times \pi_2^{-1}(s_2) \mbox{  ,  }
(s_1,s_2) \in \Scal_1 \times \Scal_2 \}$$
où $\Scal_1$ et $\Scal_2$ sont des découpages de $\Cc$ en carrés de taille
$r$. En remplaçant les courants laminaires par les courants
géométriques dans la preuve de R. Dujardin (voir Proposition 4.4 dans
\cite{Du2}), on obtient la

\begin{proposition}{\label{prop4}}

Il existe un courant $S_{\Qcal} \leq S$ uniformément géométrique dans
les $Q \in \Qcal$ tel que $M_{\Omega}(S-S_{\Qcal}) \leq Cr^2$ (où $C$ est
une constante universelle et $M_{\Omega}$ est la masse dans
$\Omega$).

\end{proposition}

Nous allons expliquer brièvement la construction de R. Dujardin car nous en
aurons besoin par la suite. On notera $S_m= \frac{1}{m} \sum_{i=0}^{m-1} \frac{[f^i(L)]}{d^i}$ et
quitte à extraire une sous-suite on pourra supposer que $S_m$ converge
vers $S$.

Dans les composantes connexes de $f^i(L)$ au-dessus d'un carré $s_j$ de
$\Scal_j$ (i.e. dans $\pi_j^{-1}(s_j) \cap f^i(L)$), on a un certain
nombre de disques qui sont des graphes au-dessus de $s_j$ (qui peuvent éventuellement s'intersecter entre eux). En moyennant les
courants d'intégration sur les graphes qui ont une aire inférieure à
$\frac{1}{2}$, on obtient ainsi deux courants $S_{m, \Qcal_j}$ inférieurs à
$S_m$ et qui vérifient:
$$\langle S_m-S_{m, \Qcal_j}, \pi_j^{*} \omega \rangle \leq Cr^2.$$
Cette inégalité résulte de la formule de Riemann-Hurwitz et du fait
que le genre de $f^i(L)$ vaut $0$.

Par ailleurs, $S_{m, \Qcal_j}$ s'écrit $\frac{1}{m} \sum_{i=0}^{m-1}
\frac{[C_{\Qcal_j}^i]}{d^i}$ où les $C_{\Qcal_j}^i$ sont des courbes à
bord dans $\pi_j^{-1}(\partial \Scal_j)$.

On réunit maintenant ces deux courants en posant:
$$S_{m, \Qcal}=\frac{1}{m} \sum_{i=0}^{m-1}
\frac{[C_{\Qcal_1}^i \cup C_{\Qcal_2}^i]}{d^i}.$$
C'est un courant uniformément géométrique dans les cubes de $\Qcal$
car les courbes $C_{\Qcal_1}^i \cup C_{\Qcal_2}^i$ qui composent $S_{m,
  \Qcal}$ sont à bord dans $\partial \Qcal$. De plus, la suite $S_{m,
  \Qcal}$ vérifie:
$$\langle S_m-S_{m, \Qcal}, \pi_1^{*} \omega + \pi_2^{*} \omega \rangle \leq Cr^2.$$
Pour les valeurs d'adhérences $S_{\Qcal}$ de $S_{m, \Qcal}$ on a
donc bien:
$$M_{\Omega}(S-S_{\Qcal}) \leq Cr^2.$$
Enfin, il est assez facile de voir que le courant $S_{\Qcal}$
est uniformément géométrique dans les $Q \in \Qcal$.

$ $

La proposition ci-dessus est le point clé pour construire une suite
$S_k$ qui croît vers $S$ telle que $T \wedge S_k$ croisse vers $T
\wedge S$. En effet, si on reprend les arguments de R. Dujardin (voir
\cite{Du2}), on constate
que pour un quadrillage $\Qcal$ de taille $r$ bien choisi (i.e. bien positionné), on a:
$$M_{\Omega}(T \wedge S - T \wedge S_{\Qcal}) \leq \epsilon(r).$$
Ici $S_{\Qcal}$ est le courant construit précédemment et
$\epsilon(r)=C \omega(G,r)^{1/3}$, où $\omega(G,r)$ est le module de
continuité du potentiel de $T$ de rayon $r$ et $C$ est une
constante universelle. En particulier $\epsilon(r)$ est une suite qui tend vers $0$ quand $r$ décroît vers $0$.

En itérant le procédé, on obtient donc le

\begin{theoreme}

Il existe une suite de subdivisions $\Qcal_k$ de $\overline{\Omega}$ et
une suite $S_k$ de courants uniformément géométriques dans les cubes
de la subdivision telles que:

i) La suite de courants $S_k$ croît vers $S$.

ii) Celle des mesures $T \wedge S_k$ croît vers $T \wedge S$.

\end{theoreme}

\section{{\bf Minoration du plus grand exposant de Lyapunov de $\nu$}}

Nous avons vu au début de cet article que l'entropie de $\nu$ est
minorée par $\log(d)$. En utilisant la formule de Ruelle
(voir le paragraphe \ref{exposant}), on obtient alors $$\int \chi^{+}
(x) d\nu(x) \geq \frac{\log(d)}{2}$$
où $\chi^{+} (x)= \sum_{i \mbox{, } \chi_i (x) > 0}  \chi_i (x) .$\\
Nous allons maintenant raffiner ce résultat. Plus précisément, en
adaptant les arguments de R. Dujardin (voir \cite{Du3}) à notre cas
qui est non inversible nous allons montrer le:
\begin{theoreme*}

Pour $\nu$ presque tout point $x$, le plus grand exposant de Lyapunov
en $x$ est supérieur ou égal à $\frac{\log(d)}{2}$.

\end{theoreme*}

Quitte à extraire une sous-suite, nous supposerons dans la suite que
$S_m=\frac{1}{m} \sum_{i=0}^{m-1} \frac{[f^i(L)]}{d^i}$ converge
vers $S$. De plus, quitte à changer la suite $\frac{1}{m}
\sum_{i=0}^{m-1} \frac{[f^i(L)]}{d^i}$ en une
suite $\frac{1}{m} \sum_{i=0}^{m-1} \frac{[f^i(L_i)]}{d^i}$ qui
converge aussi vers $S$, nous
pourrons supposer que $f^i(L)$ n'a pas de multiplicité.

Avant de donner l'idée de la preuve, remarquons que le théorème revient à démontrer que
l'ensemble $\{ x \mbox{ , } \lim_{n \rightarrow \infty} \frac{1}{n} \log \| D_x f^n  \| <
\frac{\log(d)}{2} \}$ est de mesure nulle pour $\nu$ (la limite existe par le
théorème sous-additif de Kingman). Autrement dit, il suffit de montrer
que l'ensemble
$$B_{\alpha, n_0}:= \{ x \mbox{ , } \forall n \geq n_0 \mbox{ , } \frac{1}{n} \log \| D_x f^n  \| <
(1-\alpha)\frac{\log(d)}{2} \}$$
est de mesure nulle pour $\nu$ pourvu que $\alpha$ soit petit et $n_0$ grand.

Enfin, quitte à considérer un recouvrement fini de $\Pp^2$ par des
ouverts $\Omega$ et à prendre un $m$ grand, il suffira donc de majorer
$T \wedge S_m(\Omega_n)$ par $\epsilon$ où 
$$\Omega_n:=\{ x \in \Omega \mbox{ , }  \frac{1}{n} \log \| D_x f^n  \| <
(1-\alpha)\frac{\log(d)}{2} \}.$$

$ $

Voici maintenant l'idée de la preuve.

Admettons que l'on puisse construire environ $d^{i+n}$ disques consistants
dans $f^{i+n}(L)$. Alors, si on
n'a pas de problèmes liés à l'ensemble critique, on obtient $d^{i+n}$
préimages de ces disques dans $f^i(L)$. La plupart d'entre eux ont
donc une aire inférieure à $Kd^{-n}$ (car l'aire de $f^i(L)$ vaut
$d^i$) et on a créé ainsi beaucoup de
points qui ne sont pas dans $\Omega_n$.

$ $

Le plan de la preuve du théorème sera donc le suivant: dans un premier paragraphe, nous allons construire ces disques dans
$\frac{f_{*}^n S_m}{d^n}$. Ensuite, on prendra les préimages de ces
disques par $f^n$ et nous utiliserons un argument
longueur-aire. Enfin, la
dernière partie sera consacrée à la démonstration du théorème.

\subsection{{\bf Construction des disques}}

Quitte à considérer un recouvrement fini de $\Pp^2(\Cc)$ par des
ouverts $U$, nous pouvons supposer que $f^n(\Omega_n) \subset U$.

On considère maintenant un découpage $\Qcal$ d'un voisinage de $\overline{U}$ en $4$-cubes affines de taille $\frac{1}{k}$. Dans toute la
suite le quadrillage sera considéré bien choisi par rapport à la mesure $\nu$. En particulier $\nu$ charge peu un petit voisinage du bord de $\Qcal$
(voir le lemme 4.5 dans \cite{Du2}).
Notons $\Qcal_{\lambda}=\cup Q_{\lambda}$ où les
$Q_{\lambda}$ sont les homothétiques des $Q \in \Qcal$ de
rapport $\lambda$. Quitte à remplacer $\Omega_n$ par l'ouvert $\Omega_n
\cap f^{-n}( \Qcal_{\lambda})$, nous pourrons donc supposer dans la suite que $f^n(\Omega_n) \subset
\Qcal_{\lambda}$ (cela fait changer éventuellement $m$ en un $m^{'}$
plus grand).

Maintenant, à partir du courant $\frac{f^n_{*} S_m}{d^n}=\frac{1}{m}
\sum_{i=0}^{m-1} \frac{[f^{i+n}(L)]}{d^{i+n}}$, on peut construire un
courant $R_{m, \Qcal}$ uniformément géométrique dans les cubes du
quadrillage de $\Qcal$ et qui vérifie (voir le paragraphe \ref{geometrique})
$$\left( T \wedge \frac{f^n_{*} S_m}{d^n} - T \wedge R_{m, \Qcal}
\right) (U) \leq
\epsilon(k).$$
Ici $\epsilon(k)$ ne dépend pas de $m$ et $n$. En effet d'une part on a
$$M_U \left( \frac{f^n_{*} S_m}{d^n}-R_{m, \Qcal} \right) \leq \frac{C}{k^2}$$
où $C$ est indépendante de $m$ et $n$ et d'autre part comme le quadrillage $\Qcal$
est supposé bien choisi par rapport à $\nu$, il l'est pour $T
\wedge  \frac{f^n_{*} S_m}{d^n}$ pour des $m$ grands (voir le lemme
4.5 de \cite{Du2}).

Dans la suite, on posera $R_{m, \Qcal}=\int [\Gamma] d \nu_{m, \Qcal}(\Gamma)$. Dans cette expression les $\Gamma$ sont des
disques inclus dans les cubes de $\Qcal$. Ce sont ces disques de
$R_{m, \Qcal}$ que l'on va tirer en arrière par
$f^n$.

\subsection{Tiré en arrière des disques}

Dans la suite, nous allons jeter un certain nombre de disques de
$R_{m, \Qcal}$ de sorte à pouvoir effectuer l'argument longueur-aire
et ainsi créer beaucoup de points qui ne seront pas dans $\Omega_n$. A
chaque fois que nous enlèverons un disque $\Gamma$, nous donnerons la
perte $\int_{\Gamma} T$ ainsi enregistrée.

$ $

Fixons $i$ compris entre $0$ et $m-1$.

On ne considèrera que les $\Gamma$ qui sont dans $f^{i+n}(L) \cap
Q_{\lambda}$. En particulier, comme $\Gamma$ est un graphe, $f^n(\Omega_n) \cap \Gamma$ est inclus dans un disque
$\Gamma^{'} \Subset \Gamma$ et le module de l'anneau $\Gamma
\setminus \Gamma^{'}$ est minoré par une constante indépendante de
$\Gamma$ (qui ne dépend que de $k$ et de $\epsilon$).

Comme la courbe $f^{i+n}(L)$ n'a pas de multiplicité, on peut tirer en
arrière par $f^n$ les disques $\Gamma$. On notera $(\Gamma)_j^{-n}$ et
$(\Gamma^{'})_j^{-n}$ les préimages de $\Gamma$ et $\Gamma^{'}$ qui
sont dans $f^i(L)$. Le nombre de disques $\Gamma$ pour lesquels
$(\Gamma)_j^{-n}$ est d'aire supérieure à $K d^{-n}$ est majoré par
$\frac{d^{i+n}}{K}$ (car $f^i(L)$ est d'aire $d^i$). Si on désigne par
$\Dcal(i)$ l'ensemble de ces disques et qu'on les retire à $R_{m,
  \Qcal}$, on enregistre une perte $\int_{\Dcal(i) \cap Q_{\lambda}} T
\leq C(k) \frac{d^{i+n}}{K}$.

Enfin, si on considère un disque $\Gamma$ que l'on n'a pas enlevé, on a
$$\mbox{Diamètre}((\Gamma^{'})_j^{-n}) \leq C(k, \epsilon) \sqrt{K} d^{-n/2}$$ 
grâce à un argument longueur-aire (voir l'appendice de \cite{BD2}).

Alors, quitte à considérer un disque intermédiaire entre $\Gamma^{'}$ et $\Gamma$, on obtient grâce aux inégalités de Cauchy (comme dans \cite{Du3}):
$$\| D f_{j}^{-n}(p)v(p) \| \leq C(k,K, \epsilon) d^{-n/2}$$
pour $p \in \Gamma^{'}$ et un certain vecteur unitaire $v(p)$. Ici
$f_{j}^{-n}$ désigne la branche inverse de $f^n$ qui envoie $\Gamma$
sur $(\Gamma)_j^{-n}$.

Ces points $f_{j}^{-n}(p)$ vérifient $\| D f^n(f_{j}^{-n}(p)) \| \geq
\frac{d^{n/2}}{C(k,K, \epsilon)}$ et ils ne peuvent donc pas appartenir à
$\Omega_n$ (si $n$ est suffisamment grand).

\subsection{Démonstration du théorème}

Rappelons que nous devons voir que
$$T \wedge S_m( \Omega_n ) \leq \epsilon.$$
Si $p$ est un point de $\Omega_n$, alors $f^n(p)$ est soit
dans un des disques que l'on a enlevés à $R_{m, \Qcal}$, soit dans la
partie de $\frac{f^n_{*} S_m}{d^n}$ où il n'y a pas de disque de
$R_{m, \Qcal}$. Autrement dit, en utilisant la relation $f_{*}^n(T \wedge S_m)=T \wedge \frac{f^n_{*} S_m}{d^n}$ et l'estimée de perte que l'on avait obtenue, on a:
$$T \wedge S_m(\Omega_n) \leq T \wedge \frac{f^n_{*} S_m}{d^n}
(f^n(\Omega_n)) \leq \frac{C(k)}{K}  + \epsilon(k) \leq \epsilon$$
pour $k$ puis $K$ bien choisies.

\section{{\bf Majoration du plus petit exposant de Lyapunov de $\nu$}}

Dans ce paragraphe $\nu$ désigne toujours la mesure $T \wedge S$.

Nous voulons ici montrer le

\begin{theoreme*}
Supposons que $\nu$ ne charge pas les courbes algébriques.

Alors, pour $\nu$ presque tout point $x$ hors du support de $\mu$, le plus
petit exposant de Lyapunov $\chi_1(x)$ est négatif ou nul.

\end{theoreme*}

Voici le plan de la démonstration. Dans un premier
paragraphe, nous allons montrer que la mesure de $\{ x \notin
\mbox{Support}(\mu) \mbox{ , } \chi_2(x) = \chi_1(x) > 0 \}$ est nulle. Les
exposants étant égaux, nous pourrons en particulier utiliser les travaux de F. Berteloot et C. Dupont (voir
\cite{BDu}) sur les linéarisations le long d'orbites négatives. Dans
le second paragraphe on verra que la mesure de $\{ x \notin
\mbox{Support}(\mu) \mbox{ , } \chi_2(x) > \chi_1(x) > 0 \}$ vaut $0$. Cette fois-ci les exposants étant disjoints, nous aurons une
direction stable (donnée par le plus petit exposant) et une direction
instable. Dans ce contexte, nous pourrons en particulier utiliser la
tranformée de graphe.

\subsection{{\bf Majoration de la mesure de  $ \{ x \not \in
    \mbox{Support} (\mu) \mbox{ , } \chi_2 (x) = \chi_1 (x) > 0 \} $ }}

Avant de passer à l'estimation de cette mesure, nous allons
faire quelques rappels sur le procédé de branches inverses de
J.-Y. Briend et J. Duval (voir \cite{BD1}) ainsi que sur la linéarisation de
F. Berteloot et C. Dupont (voir \cite{BDu}).

\subsubsection{Rappels}{\label{inverse}}

On reprend les notations du paragraphe \ref{exposant}. En particulier
$\tau_x$ désigne une bonne famille de cartes holomorphes de $\Pp^2$ et
$f_x= \tau_{f(x)}^{-1} \circ f \circ \tau_x$. Dans la suite nous
noterons aussi:
$$f_x^n=\tau_{f^n(x)}^{-1} \circ f^n \circ \tau_x=f_{f^{n-1}(x)} \circ
\cdots \circ f_x$$
et pour $\widehat{x}$ dans l'extension naturelle $\widehat{\Pp^2(\Cc)}$,
$$f_{\widehat{x}}^{-n}= f_{x_{-n}}^{-1} \circ \cdots \circ
f_{x_{-1}}^{-1}$$
(quand cette expression est bien définie). On notera aussi
$\widehat{X}$ l'ensemble des bons points de Pesin pour $\widehat{\nu}$
(voir le paragraphe \ref{exposant}).

Le borélien $\widehat{X}^{'}= \{ \widehat{x} \in \widehat{X} \mbox{ , }
\chi_1(\widehat{x})= \chi_2(\widehat{x}) \mbox{ , }
\chi_1(\widehat{x}) > 0 \} $ est invariant par
$\widehat{f}$. En particulier, si on considère la mesure
$\widehat{\nu^{'}}$ définie par $\widehat{\nu^{'}}(A)=\widehat{\nu}(A
\cap \widehat{X}^{'})$, on obtient une mesure invariante pour laquelle
$\log ^{+} \| (D \widehat{f})^{-1} \|$ est intégrable (voir
l'appendice $A.1$ de \cite{K}).

En particulier nous pouvons utiliser d'une part les résultats de J.-Y. Briend et J. Duval
(voir \cite{BD1} et \cite{B}) et d'autre part le procédé de linéarisation de
F. Berteloot et C. Dupont (voir \cite{BDu}) pour la mesure
$\widehat{\nu^{'}}$ et ainsi obtenir le

\begin{theoreme*}{\cite{BD1} et \cite{BDu}}.

Soient $\gamma$ et $R$ suffisamment petits.
Il existe un borélien $\widehat{Y}^{'} \subset \widehat{X}^{'}$ avec
$\widehat{\nu}(\widehat{X}^{'} - \widehat{Y}^{'})=0$ et des fonctions
mesurables $\eta: \widehat{Y}^{'} \rightarrow ]0,R]$,
    $F: \widehat{Y}^{'} \rightarrow ]0, + \infty[$,
        $C : \widehat{Y}^{'} \rightarrow [1, + \infty[$
            et $S: \widehat{Y}^{'} \rightarrow ]0,R]$
        tels que pour tout $\widehat{x} \in \widehat{Y}^{'}$ et
        tout $n$ dans $\Nn$ on
        ait:

1) $f_{\widehat{x}}^{-n}$ est définie sur $B(0, \eta(\widehat{x}))$.

2) $\mbox{Lip} f_{\widehat{x}}^{-n} \leq C(\widehat{x})
e^{-n(\chi_1(\widehat{x})- \gamma/2)}$.

3) $\| D_0 f^n_{x_{-n}} \| \leq F(\widehat{x})
   e^{n(\chi_2(\widehat{x}) + \gamma)}$.

4) $S \leq \eta$ et $D_0 f_{\widehat{x}}^{-n} (B(0,S(\widehat{x})))
\subset f_{\widehat{x}}^{-n} B(0, \eta(\widehat{x}))$.

\end{theoreme*}

\begin{remarque*}

Les deux premiers points proviennent du fait que $\chi_1(\widehat{x}) >
0$ (voir \cite{BD1} et \cite{B}). Le troisième est une conséquence
directe de la théorie de Pesin. Enfin le dernier point est le seul qui
utilise que $ \chi_1(\widehat{x}) = \chi_2(\widehat{x})$ (voir
\cite{BDu}).

\end{remarque*}

\subsubsection{Majoration de la mesure de $\{ x \not \in
\mbox{Support}(\mu) \mbox{ , } \chi_2(x) = \chi_1(x) > 0 \}$}

Voici le plan de cette majoration.

Après un premier paragraphe consacré à des préliminaires, nous verrons
dans le second que montrer que la mesure de $E_1=\{ x \notin
\mbox{Support}(\mu) \mbox{ , } \chi_2(x) = \chi_1(x) > 0 \}$ est nulle revient à calculer la
mesure pour $\nu$ de préimages $f_i^{-n}(B)$. Enfin dans le troisième
paragraphe nous effectuerons ce calcul.

\subparagraph{1) Préliminaires}

$\mbox{ }$

Commençons par ramener la majoration de $E_1$ à celle d'un borélien
$\widehat{E}_1(n)$ qui aura de bonnes propriétés de récurrence et
d'uniformité.

Remarquons tout d'abord qu'il nous suffit de voir que $\nu( \{ x \in
E_1 \mbox{ , } \chi_1(x) \geq \beta > 0  \})$ est nul (avec $\beta$ petit). On notera toujours $E_1$ cet ensemble. Par ailleurs,
$$\nu(E_1)= \widehat{\nu}(\pi^{-1}(E_1))= \widehat{\nu}(\pi^{-1}(E_1)
\cap \widehat{X})$$
où $\pi$ est la projection de $\widehat{\Pp^2(\Cc)}$ sur $\Pp^2(\Cc)$
et $\widehat{X}$ est l'ensemble des bons points de Pesin (voir le
paragraphe \ref{exposant}).

Ce dernier ensemble $\pi^{-1}(E_1) \cap \widehat{X}$ est inclus dans
$$\widehat{E}_1=\{ \widehat{x} \in \widehat{X} \mbox{ , } \pi(\widehat{x}) \notin
\mbox{Support}(\mu) \mbox{ , } \chi_2(\widehat{x}) =
\chi_1(\widehat{x}) \geq \beta \}$$
(car la définition des exposants de Lyapunov ne dépend que des
trajectoires positives).

La majoration de $\nu(E_1)$ se ramène donc à celle de
$\widehat{\nu}(\widehat{E}_1)$.

Ensuite, si on reprend les notations du théorème précédent, on constate
que quitte à remplacer $\widehat{E}_1$ par 
$$\widehat{E}_{1,p}= \{
\widehat{x} \in \widehat{E}_1 \cap \widehat{Y}^{'} \mbox{ , }
d(\pi(\widehat{x}),\mbox{Support}(\mu)) \geq \frac{1}{p}   \mbox{ , }     S(\widehat{x}) \geq \frac{1}{p} \mbox{ , } C(\widehat{x})
\leq p \mbox{ , } F(\widehat{x}) \leq p \}$$
pour $p$ grand, on pourra
supposer dans la suite que nous avons ces contrôles.

Enfin, l'ensemble $\widehat{E}_1$ aura de bonnes
propriétés de récurrence en utilisant le

\begin{lemme}{\label{recurrence}}

$$(\widehat{\nu}(\widehat{E}_1))^2 \leq \lim_{m \rightarrow  +  \infty}
  \frac{1}{m} \sum_{n=0}^{m-1} \widehat{\nu}(
  \widehat{f}^n(\widehat{E}_1) \cap \widehat{E}_1).$$

\end{lemme}

\begin{proof}

Il suffit d'utiliser d'une part la décomposition de $\widehat{\nu}$ en mesures
ergodiques $\widehat{\nu}_{\alpha}$ et d'autre part le fait qu'une mesure ergodique $\widehat{\nu}_{\alpha}$
vérifie 
$$\lim_{m \rightarrow +  \infty} \frac{1}{m} \sum_{n=0}^{m-1}
\widehat{\nu}_{\alpha}(\widehat{f}^n(\widehat{E}_1) \cap
      \widehat{E}_1)=(\widehat{\nu}_{\alpha}(\widehat{E}_1))^2.$$ 

\end{proof}

En effet, grâce à ce lemme, pour démontrer que
 $\widehat{\nu}(\widehat{E}_1)$ (et donc $\nu(E_1)$) est de mesure
 nulle, il nous suffira de montrer que
 $\widehat{\nu}(\widehat{E}_1(n))$ est petit pour $n$ assez grand (avec $\widehat{E}_1(n)=\widehat{f}^n(\widehat{E}_1) \cap
      \widehat{E}_1$). $\widehat{E}_1(n)$ est le bon borélien que l'on
 cherchait.

Maintenant, grâce au procédé de linéarisation (point $4$ du théorème
précédent), nous avons la

\begin{proposition*}
Il existe $\rho > 0$ tel que:

$$\exists n_0 \mbox{ , } \forall n \geq n_0 \mbox{ , }  \forall
\widehat{x} \in \widehat{E}_1$$ on a:
$$\mbox{Diamètre intérieur de } f_{x_{-n}}^n( D_0
f_{\widehat{x}}^{-n}(B(0,\frac{1}{4p}))) \geq \rho.$$

\end{proposition*}

\begin{proof}

Faisons un raisonnement par l'absurde.

Si la proposition est fausse, on obtient des suites $n_l$ et $\widehat{x}_l \in  \widehat{E}_1$ telles
que le diamètre intérieur de $f_{(x_l)_{-n_l}}^{n_l}( D_0
f_{\widehat{x}_l}^{-n_l}(B(0,\frac{1}{4p})))$ soit majoré par
$\frac{1}{l}$.

Par ailleurs la fonction $g_{n_l}(u)=f_{(x_l)_{-n_l}}^{n_l} D_0
f_{\widehat{x}_l}^{-n_l}(u)$ est définie sur $B(0,\frac{1}{p})$ et son
image est incluse dans $B(0,R)$ (voir le point $4$ du théorème
précédent). Quitte à extraire une sous-suite, $g_{n_l}$ converge donc
vers une fonction holomorphe $g: B(0,\frac{1}{p}) \mapsto
\overline{B(0,R)}$.

Maintenant, on a d'une part $D_0 g$ qui est égal à $I$ et d'autre part le diamètre
intérieur de  $g(B(0,\frac{1}{4p}))$ qui vaut $0$. On obtient alors
une contradiction par le théorème d'inversion locale.

\end{proof}

\subparagraph{2) Branches inverses}

$\mbox{ }$

Commençons par recouvrir $\Pp^2$ par des boules $B$ telles que
$\tau_x^{-1}(B) \subset B(0, \rho)$ pour tout $x$ de $B$ ($\rho$ est
une constante petite qui vérifie les conclusions de la proposition
précédente).

Ensuite nous appellerons bonne composante de $f^{-n}(B)$ un ouvert de
la forme $\tau_{y_{-n}} \circ f_{\widehat{y}}^{-n} \circ
\tau_{y_0}^{-1}(B)$ avec $y_0 \in B$ et $\widehat{y} \in
\widehat{E}_1(n)$.

Nous avons alors (comme dans \cite{BD1}):
$$\widehat{\nu} ( \widehat{E}_1(n)
\cap \pi^{-1}(B)) = \widehat{\nu} (\widehat{f}^{-n}( \widehat{E}_1(n)
\cap \pi^{-1}(B))) \leq \nu( \pi (\widehat{f}^{-n}( \widehat{E}_1(n)
\cap \pi^{-1}(B))))$$
car $\pi_{*} \widehat{\nu} = \nu$.

Mais $\pi (\widehat{f}^{-n}( \widehat{E}_1(n)
\cap \pi^{-1}(B)))$ est inclus dans l'union des bonnes composantes de
$f^{-n}(B)$ (que nous noterons $\cup \mbox{BC}$) c'est-à-dire
$$ \widehat{\nu} ( \widehat{E}_1(n)
\cap \pi^{-1}(B)) \leq  \nu ( \cup \mbox{BC}) .$$
Par ailleurs, remarquons qu'une bonne composante ci-dessus évite un
$\frac{1}{2p}$-voisinage du support de $\mu$ (car d'une part les bonnes
composantes contiennent un point qui est hors d'un
$\frac{1}{p}$-voisinage du support de $\mu$ par définition de
$\widehat{E}_1(n)$ et d'autre part les bonnes composantes sont de diamètre
exponentiellement petit par le point $2$ du théorème du paragraphe \ref{inverse}).

\subparagraph{3) Calcul de $\nu(\cup \mbox{BC})$}

$\mbox{ }$

Cette opération va se dérouler en trois étapes.

Dans la première, nous allons faire des simplifications géométriques:
nous remplacerons $\nu$ par $T \wedge S_{\Qcal}$ où $S_{\Qcal}$ est un
courant uniformément géométrique et on verra que grâce à la
linéarisation, les bonnes composantes seront presque des
ellipsoïdes.

Après cette opération, nous aurons à démontrer trois faits. Le premier
est que $T \wedge [ \Gamma] (\mbox{BC})      \leq Cd^{-n}$ pour un disque
$\Gamma$ de $S_{\Qcal}$ et $\mbox{BC}$ une bonne composante. Le second
est que la mesure des graphes (pour la mesure sur les graphes de $S_{\Qcal}$) qui coupent une bonne
composante tend vers $0$ quand $n$ croît vers l'infini. Enfin, nous
verrons que le nombre de bonnes composantes de $f^{-n}(B)$ qui évitent le
petit voisinage du support de $\mu$ est majoré par $Cd^n$.

Après avoir démontré ces trois faits, nous aurons que $T \wedge S_{\Qcal}(\cup \mbox{BC})$ sera majoré par
$Cd^n T \wedge S_{\Qcal}(\mbox{BC}) \leq Cd^n Cd^{-n} \epsilon(n)$ qui
est la majoration cherchée.

\subparagraph{i) Simplifications géométriques}

Quitte à considérer un recouvrement fini de $\Pp^2(\Cc)$ par des
ouverts $\Omega$, nous pouvons supposer que $\cup \mbox{BC}$ est
inclus dans $\Omega$. On peut construire un
découpage d'un voisinage de $\overline{\Omega}$ en $4$-cubes affines
et un courant $S_{\Qcal}= \int [\Gamma] d \nu_{\Qcal}(\Gamma)$ uniformément géométrique dans les cubes de
ce quadrillage tels que la différence $(T \wedge S - T \wedge
S_{\Qcal})(\Omega)$ soit aussi petite que l'on veut (voir le
paragraphe \ref{courant}). Nous pouvons supposer de plus que $\nu$
charge très peu un petit voisinage du bord de $\Qcal$ et alors notre
majoration de $\nu(\cup \mbox{BC})$ se ramène à celle de $T \wedge
S_{\Qcal}(\cup \mbox{BC})$ pour les bonnes composantes qui évitent un
petit voisinage du bord de $\Qcal$. Comme précédemment, nous noterons
toujours $\cup \mbox{BC}$ l'union de ces bonnes composantes.

Expliquons maintenant pourquoi les bonnes composantes ressemblent à des ellipsoïdes. Notons ici $B_i^{-n}= \tau_{y_{-n}} \circ
f_{\widehat{y}}^{-n} \circ \tau_{y_0}^{-1}(B)$ une des bonnes
composantes de $f^{-n}(B)$.

D'après la proposition précédente, nous savons que
$$\tau_{y_0}^{-1}(B) \subset B(0,\rho) \subset f_{y_{-n}}^n D_0
f_{\widehat{y}}^{-n} B(0, \frac{1}{4p}).$$
On en déduit donc que $B_i^{-n} \subset \tau_{y_{-n}} D_0
f_{\widehat{y}}^{-n} B(0, \frac{1}{4p})$ et $D_0
f_{\widehat{y}}^{-n} B(0, \frac{1}{4p})$ est bien un ellipsoïde.

\subparagraph{ii) Démonstration des trois faits}

\subparagraph{$\alpha$) Calcul de $T \wedge [ \Gamma] (\tau_{y_{-n}} D_0
f_{\widehat{y}}^{-n} B(0, \frac{1}{4p}))$}

Par construction $\Gamma$ est un graphe (i.e. de la forme
$(z,\phi(z))$) et le diamètre de $ D_0
f_{\widehat{y}}^{-n} B(0, \frac{1}{p})$ est majoré par $e^{-
  \chi_1(\widehat{y})n + \gamma n}$.

Par ailleurs, nous avons des contrôles uniformes d'une part sur les
$\tau_{y_{-n}}^{-1}$ et d'autre part sur les dérivées de $\phi$ (via
le théorème de Cauchy car $\tau_{y_{-n}} D_0
f_{\widehat{y}}^{-n} B(0, \frac{1}{p})$ est un peu loin du bord de
$\Qcal$). La partie $\Gamma_0$ de $\tau_{y_{-n}}^{-1}(\Gamma)$ susceptible
de rencontrer $D_0 f_{\widehat{y}}^{-n} B(0, \frac{1}{p})$ se met
donc sous la forme d'un graphe $(z,\psi(z))$ où $z$ est dans un disque
$D_0$ de taille $e^{-  \chi_1(\widehat{y})n + \gamma n}$ et
$|\psi^{'}(z)| \leq e^{-  \chi_1(\widehat{y})n + \gamma n}$ sur
$D_0$ (on s'est placé dans un repère où l'axe des abscisses est tangent à $\Gamma_0$ en un point de $ D_0 f_{\widehat{y}}^{-n}
B(0, \frac{1}{4p})$).

Maintenant le calcul de $\int_{\tau_{y_{-n}} D_0
f_{\widehat{y}}^{-n} B(0, \frac{1}{4p})} T \wedge [ \Gamma]=\int_{ D_0
f_{\widehat{y}}^{-n} B(0, \frac{1}{4p})} \tau_{y_{-n}}^{*} T \wedge
[ \Gamma_0]$ va se faire par un argument de capacité.

En effet, notons $\Ecal_i=D_0
f_{\widehat{y}}^{-n} B(0, \frac{1}{ip})$ ($i=1,2,3,4$),
$D^{'}=\pi(\Gamma_0 \cap \Ecal_4)$ et $D=\pi(\Gamma_0 \cap \Ecal_1)$
(où $\pi$ est la projection du graphe $\Gamma_0$ sur l'axe des
abscisses).

Alors, si on a une fonction $G$ psh dans $\tau_{y_{-n}}(\Ecal_1)$ avec
$0 \leq G \leq Cd^{-n}$ dans cet ensemble et $dd^cG=T$, on aura:
$$\int_{\Ecal_4} \tau_{y_{-n}}^{*} T \wedge [ \Gamma_0]= \int
\tau_{y_{-n}}^{*} T \wedge [ \pi^{*} D^{'}]= \int_{D^{'}} dd^c G \circ
\tau_{y_{-n}} \circ \pi^{-1}$$
(on considère $\pi$ comme un biholomorphisme de $\Gamma_0 \cap
\Ecal_4$ sur $D^{'}$).

Enfin comme $0 \leq G \leq C d^{-n}$ sur $\tau_{y_{-n}}(\Ecal_1)$ on a $0
\leq G \circ \tau_{y_{-n}} \circ \pi^{-1} \leq C d^{-n}$ sur $D$ et
alors
$$\int_{\tau_{y_{-n}}(\Ecal_4)} T \wedge [ \Gamma] \leq
\mbox{Cap}(D^{'},D) C d^{-n}.$$
Pour obtenir la majoration de $\int_{\tau_{y_{-n}}(\Ecal_4)} T \wedge
[ \Gamma]$ par $Cd^{-n}$ il reste donc d'une part à voir que
$\mbox{Cap}(D^{'},D)$ est bornée par une constante uniforme et d'autre
part à construire la fonction psh $G$ qui est comprise entre $0$ et
$Cd^{-n}$ sur $\tau_{y_{-n}}(\Ecal_1)$ et qui vérifie $dd^c G=T$.

$ $

Majoration de $\mbox{Cap}(D^{'},D)$:

Notons $\Delta$ l'axe des abscisses qui est tangent à $\Gamma_0$. On a
d'une part $D^{'} \subset \Delta \cap \Ecal_3$ et d'autre part $\Delta
\cap \Ecal_2 \subset D$ (car $| \psi^{'} | \leq e^{-
  \chi_1(\widehat{y})n + \gamma n}$, $\| D_0 f_{y_{-n}}^n \| \leq p e^{
  \chi_2(\widehat{y})n + \gamma n}$ et $2 \chi_1(\widehat{y}) >
\chi_2(\widehat{y})$).

Cela implique donc que $\mbox{Cap}(D^{'},D) \leq \mbox{Cap}(\Delta
\cap \Ecal_3, \Delta \cap \Ecal_2)$ et il est clair que cette
dernière quantité est majorée par une constante uniforme.

$ $

Construction de $G$:

Notons $\pi_0$ la projection de $\Cc^3$ sur $\Pp^2$ et $G_0 \circ
\sigma$ la fonction de Green de $f$ définie sur une boule $\tau_{y_0} B(0,\eta(\widehat{y}))$ qui contient $B$ (ici $\sigma$ est une
section locale de $\pi_0$).

$G_0 \circ \sigma$ vérifie $dd^c (G_0 \circ \sigma)=T$ et $| G_0 \circ
\sigma | \leq \frac{C}{2}$ (voir \cite{FS1}). Par ailleurs on notera $f_i^{-n}$ la branche inverse de $f^n$ telle que $\tau_{y_{-n}} \Ecal_1
\subset f_i^{-n}(\tau_{y_0}B(0,\eta(\widehat{y})))$ (voir le point $4$
du théorème de la section \ref{inverse}).

Maintenant, en utilisant le lemme 2.3 de Ueda (voir \cite{Ue}), on
peut construire une section $\sigma_n$ de $\pi_0$ telle que $F^n \circ
\sigma_n = \sigma \circ f^n$ sur $f_i^{-n}(\tau_{y_0}B(0,\eta(\widehat{y})))$ (ici $F$ est un relevé
polynomial de $f$ par $\pi_0$). La fonction $G= G_0 \circ \sigma_n +
\frac{C}{2d^n}$ vérifie alors les conditions que l'on voulait.

\subparagraph{$\beta$) Majoration de la mesure des graphes de
  $S_{\Qcal}$ qui coupent $\tau_{y_{-n}}(\Ecal_4)$}

Il s'agit ici de voir que $\nu_{\Qcal}( \{ \Gamma \mbox{ qui coupent }
\tau_{y_{-n}}(\Ecal_4) \} ) \leq \alpha$ pour $n$ grand.

Remarquons tout d'abord que nous pouvons enlever les bonnes
composantes qui sont dans un petit voisinage de $\{ \nu(S,x) \geq
\alpha \}$ (où $\nu(S,x)$ est le nombre de Lelong de $S$ en $x$). En
effet, $\nu$ ne chargeant pas  $\{ \nu(S,x) \geq \alpha \}$, ce
voisinage $V$ est de mesure aussi petite que l'on veut. Il s'agira
alors dans la suite de majorer $T \wedge S_{\Qcal} (\cup \mbox{BC})$
où $\cup \mbox{BC}$ désigne les bonnes composantes qui évitent $V$.

Maintenant, si on considère un point $x$ hors de $V$, on a l'existence
d'un $r_0(x)$ pour lequel $S \wedge \omega (B(x,r)) < \alpha r^2$ si
$r \leq r_0(x)$ (nous supposerons de plus que $r_0(x)$ est très petit
devant le diamètre intérieur de $V$ et la distance minimale entre une
bonne composante et le bord de $\Qcal$).

En utilisant un argument de compacité sur $\Pp^2 \setminus V$, on en
déduit que toute bonne composante est incluse dans une boule $B(x,
\frac{r_0(x)}{2})$ pourvu que $n$ soit grand. On a donc bien
$\nu_{\Qcal} (\Gamma \mbox{ qui coupent }
\tau_{y_{-n}}(\Ecal_4)) < \frac{4 \alpha}{\pi}$ car dans le cas
contraire $S_{\Qcal} \wedge \omega (B(x,r_0(x))) \geq \frac{4
  \alpha}{\pi} \frac{r_0(x)^2}{4} \pi=r_0(x)^2 \alpha$ par le
théorème de Lelong.
 
\subparagraph{$\gamma$) Majoration du nombre de bonnes
  composantes}{\label{nbe}}

Si on fixe une boule $B \subset \tau_x(B(0, \rho))$ ($x \in B$), on veut majorer le nombre $N$ de bonnes
composantes de $f^{-n}(B)$ hors d'un $\frac{1}{2p}$-voisinage du
support de $\mu$ par $Cd^n$.

Remarquons que tout point $a$ de $B$ a au moins $N$ antécédents hors
du $\frac{1}{2p}$-voisinage du support de $\mu$.

Nous allons maintenant majorer ce nombre d'antécédents en utilisant un
résultat de J.E. Forn{\ae}ss et N. Sibony (voir \cite{FS2}).

Si $\Phi$ est une fonction $C^{\infty}$ à support compact, comprise
entre $0$ et $1$, et qui vaut $1$ sur un $\frac{1}{3p}$-voisinage du
support de $\mu$ et $0$ hors d'un $\frac{1}{2p}$-voisinage, on a (voir
le lemme 6.4 de \cite{FS2})
$$\mbox{Cap} \left( \left\{ a \mbox{ , } |\lambda_n(a, \Phi)-\int \Phi d \mu | \geq s \right\} \cap
B(0, \frac{1}{2}), B(0,1) \right) \leq \frac{M |\Phi|_2}{s d^n}.$$
Ici $B(0, \frac{1}{2})$ et $B(0,1)$ sont centrées en $x$ et
$\lambda_n(a, \Phi)= \frac{1}{d^{2n}} \sum_{f^n(a_i)=a} \Phi(a_i)$.

En particulier, si $a$ est un point de $B$ on a $\lambda_n(a, \Phi) \leq \frac{1}{d^{2n}}
(d^{2n} - N)$ d'où 
$$| \lambda_n(a, \Phi) - \int \Phi d \mu | = |\lambda_n(a, \Phi) - 1|
\geq \frac{N}{d^{2n}}.$$
En conclusion, nous obtenons (si le rayon de $B$ est supérieur à $\frac{\rho}{C}$):

$$\frac{d^{2n}M |\Phi|_2}{N d^n} \geq \mbox{Cap} \left( \left\{ a \mbox{ , } |\lambda_n(a, \Phi)-\int \Phi d \mu | \geq \frac{N}{d^{2n}} \right\} \cap
B(0, \frac{1}{2}), B(0,1) \right) \geq \mbox{Cap}(B,B(0,1)),$$
d'où
$$\frac{d^{2n}M |\Phi|_2}{N d^n} \geq \left( \frac{1}{\log( \rho / C)} \right) ^2,$$
qui conduit à la majoration de $N$ par la quantité souhaitée.

\subparagraph{iii) Majoration de la mesure de $E_1$}

Il s'agit ici de majorer $T \wedge S_{\Qcal} (\cup \mbox{BC})$.

Mais:
$$T \wedge S_{\Qcal} (\cup \mbox{BC}) \leq  Cd^n \max_{\mbox{BC}} T \wedge S_{\Qcal} ( \mbox{BC}),$$
car le nombre de bonnes composantes est majoré par $Cd^n$, d'où
$$T \wedge S_{\Qcal} (\cup \mbox{BC}) \leq Cd^n Cd^{-n}
\max_{\mbox{BC}} (\nu_{\Qcal}(\{ \Gamma \mbox{ qui coupent }
\mbox{BC} \})),$$
car $T \wedge [\Gamma](\mbox{BC}) \leq Cd^{-n}$, ce qui implique bien que
$$T \wedge S_{\Qcal} (\cup \mbox{BC}) \leq \epsilon$$
dès que $n$ est grand puisque $\max_{\mbox{BC}} (\nu_{\Qcal}(\{ \Gamma \mbox{ qui coupent }
\mbox{BC} \} ))$ est aussi petit que l'on veut.

\subsection{{\bf Majoration de la mesure de $\{ x \not \in
\mbox{Support}(\mu) \mbox{ , } \chi_2(x) > \chi_1(x) > 0 \}$}}

Avant de démontrer que la mesure de cet ensemble est nulle, nous allons commencer par quelques préliminaires.

\subsubsection{{\bf Préliminaires}}{\label{a}}

Le but de ce paragraphe est de ramener la majoration de la mesure de
$E_2=\{ x \notin
\mbox{Support}(\mu) \mbox{ , } \chi_2(x) > \chi_1(x) > 0 \}$ à celle de $E_2(i,j)$ où $E_2(i,j)$ est un ensemble qui
aura de bonnes propriétés d'uniformité et de récurrence.

Remarquons tout d'abord qu'il nous suffit de voir que $\nu( \{ x \in
E_2 \mbox{ , } \chi_1(x) \geq \beta > 0 \mbox{ , } \chi_2(x) - \chi_1(x)
\geq \beta \})$ est nul (avec $\beta$ petit). On notera toujours $E_2$ cet ensemble. Par ailleurs,
$$\nu(E_2)= \widehat{\nu}(\pi^{-1}(E_2))= \widehat{\nu}(\pi^{-1}(E_2)
\cap \widehat{X})$$
où $\pi$ est la projection de $\widehat{\Pp^2(\Cc)}$ sur $\Pp^2(\Cc)$
et $\widehat{X}$ est l'ensemble des bons points de Pesin (i.e. les
points qui vérifient les conclusions du premier théorème du paragraphe
\ref{exposant}).

Ce dernier ensemble $\pi^{-1}(E_2) \cap \widehat{X}$ est inclus dans
$$\widehat{E}_2=\{ \widehat{x} \in \widehat{X} \mbox{ , } \pi(\widehat{x}) \notin
\mbox{Support}(\mu) \mbox{ , } \chi_1(\widehat{x}) \geq \beta > 0 \mbox{ ,
} \chi_2(\widehat{x}) - \chi_1(\widehat{x}) \geq \beta \}$$
(car la définition des exposants de Lyapunov ne dépend que des
trajectoires positives).

La majoration de $\nu(E_2)$ se ramène donc à celle de
$\widehat{\nu}(\widehat{E}_2)$.

Ensuite, puisque sur $\widehat{E}_2$ les exposants sont strictement
positifs, on peut construire $\widehat{Y}^{'} $ avec $\widehat{\nu}
(\widehat{E}_2 - \widehat{Y}^{'})=0$ pour lequel tous les points $\widehat{x}
\in \widehat{Y}^{'} $ vérifient les points $1$ et $2$ du théorème du
paragraphe \ref{inverse} ainsi que les conclusions du deuxième
théorème du paragraphe \ref{exposant}.

En reprenant alors les notations de ce théorème ainsi que celles du
paragraphe \ref{exposant}, on constate que quitte à remplacer $\widehat{E}_2$ par 
$$\widehat{E}_{2,p}= \{
\widehat{x} \in \widehat{E}_2 \cap \widehat{Y}^{'} \mbox{, }
d(\pi(\widehat{x}),\mbox{Support}(\mu)) \geq \frac{1}{p}   \mbox{, }      C(\widehat{x}) \leq p
\mbox{, } \eta(\widehat{x}) \geq  \frac{1}{p} \mbox{, }  \delta(\widehat{x}) \mbox{ et }   \| C_{\gamma}^{\pm 1}
(\widehat{x}) \| \in [\frac{1}{p},p]  \}$$
pour $p$ grand, on pourra supposer dans la suite que nous avons ces
contrôles.

$ $

Nous allons maintenant donner des propriétés de récurrence à
$\widehat{E}_2$.

En utilisant le lemme \ref{recurrence} de la partie précédente en
temps négatif au lieu du temps positif, on
en déduit que montrer que $\widehat{\nu}(\widehat{E}_2)$ est nul découle du fait que $\widehat{\nu}(\widehat{f}^{-i}(\widehat{E}_2)
\cap \widehat{E}_2)$ est petit pour $i$ assez grand. De même, si on note $\widehat{E}_2(i)=\widehat{f}^{-i}(\widehat{E}_2)
\cap \widehat{E}_2$ et que l'on utilise à nouveau ce lemme en temps négatif, nous sommes ramenés à montrer que
$\widehat{\nu}(\widehat{f}^{-j}(\widehat{E}_2(i)) \cap
\widehat{E}_2(i))$ est petit avec $i$, $j$ assez grands, pour démontrer que $\nu(E_2)$ est nul.

Enfin, si on projette sur $\Pp^2$ on obtient:
$$\widehat{\nu}(\widehat{f}^{-j}(\widehat{E}_2(i)) \cap
\widehat{E}_2(i)) \leq \nu(\pi(\widehat{f}^{-j}(\widehat{E}_2(i)) \cap
\widehat{E}_2(i)))= \nu(E_2(i,j))$$
où
$$E_2(i,j)=\{ x_0 \mbox{ , } \exists \widehat{x} \in \widehat{f}^{-j}(\widehat{E}_2(i)) \cap
\widehat{E}_2(i) \mbox{ avec } \pi(\widehat{x})=x_0 \}.$$

\begin{remarque*}

Si on veut éviter les problèmes de mesurabilité de $E_2(i,j)$, on
pourra remplacer dans la suite $E_2(i,j)$ par son adhérence.

\end{remarque*}

L'ensemble $E_2(i,j)$ est celui que l'on cherchait et par ce que l'on
a fait, montrer que $\nu(E_2)$ est nul découle du fait que
$\nu(E_2(i,j))$ est petit.

$ $

Avant de passer à la majoration de $\nu(E_2(i,j))$ faisons une
dernière modification.

Si $x$ est un point de $\pi(\widehat{E}_2)$, on notera $E_S(x)$ l'image de $E_1(\widehat{x})$ par $D_0
\tau_x$ (ici $\widehat{x} \in \widehat{E}_2$ et
$\pi(\widehat{x})=x$). L'espace $E_S(x)$ est la direction stable en
$x$. Par ailleurs cet espace est bien défini car si
$x=\pi(\widehat{x})=\pi(\widehat{y})$ avec $\widehat{x} \mbox{ , }
\widehat{y} \in
\widehat{E}_2$, alors $E_1(\widehat{x})=E_1(\widehat{y})$ (en effet en
$\widehat{x}$ et $\widehat{y}$, les deux exposants de Lyapunov
sont distincts).

Par le théorème de Lusin, il existe un compact $K$ de mesure presque
pleine pour $\nu$ sur lequel $x \mapsto E_S(x)$ est continue. Alors
dans la suite, on cherchera à majorer la mesure de $E_2(i,j) \cap
f^{-j}(K)$ pour $\nu$. Nous noterons toujours $E_2(i,j)$ cet
ensemble.

\subsubsection{{\bf Majoration de $\nu(E_2(i,j))$}}

Avant de donner l'idée de la preuve, nous allons faire la réduction
géométrique habituelle. Quitte à considérer un recouvrement fini de
$\Pp^2(\Cc)$ par des ouverts $\Omega$, on peut supposer que $E_2(i,j)
\subset \Omega$. Maintenant, on peut construire un découpage $\Qcal$ d'un voisinage de $\overline{\Omega}$ en $4$-cubes
affines et un courant $S_{\Qcal}$ uniformément géométrique dans les
cubes de $\Qcal$ tels que la différence $(T \wedge S- T \wedge
S_{\Qcal})(\Omega)$ soit aussi petite que l'on veut. Par ailleurs, si
$\Qcal$ est bien choisi, $\nu$ charge peu un petit voisisnage du bord
de $\Qcal$. Autrement dit, si on note $Q$ un cube de $\Qcal$ et
$Q_{\lambda}$ son homothétique de rapport $\lambda$, il suffira de
majorer dans la suite
$$T \wedge S_{\Qcal}(E_2(i,j) \cap Q_{\lambda})= \int
\int_{E_2(i,j) \cap Q_{\lambda}} T \wedge [\Gamma] d \nu_{\Qcal}
(\Gamma)$$
$$= \int \int_{E_2(i,j) \cap Q_{\lambda}} \frac{f^{j*}T}{d^j} \wedge [\Gamma] d \nu_{\Qcal}
(\Gamma)= \frac{1}{d^j} \int \int_{f^j(E_2(i,j) \cap Q_{\lambda} \cap \Gamma)} T \wedge [f^j(\Gamma)] d \nu_{\Qcal}
(\Gamma).$$

\begin{remarque*}

La dernière écriture est un peu abusive car on pourrait avoir de la
multiplicité. Cependant, comme celle-ci ne joue aucun rôle, on fera cet abus: les courbes $f^j(\Gamma)$ seront supposées
sans multiplicité dans toute la suite de la démonstration.

\end{remarque*}

Voici maintenant le plan de la preuve.

Dans un premier paragraphe, nous allons construire des disques dans
les courbes $f^j(\Gamma)$. Ainsi, nous pourrons approcher le courant
$\frac{[f^j(\Gamma)]}{d^j}$ par un courant $\frac{1}{d^j} \sum
[D_{\alpha}(\Gamma)]$ uniformément géométrique dans un quadrillage que
l'on précisera. La majoration de 
$$\frac{1}{d^j} \int
\int_{f^j(E_2(i,j) \cap Q_{\lambda} \cap \Gamma)} T \wedge [f^j(\Gamma)] d \nu_{\Qcal}
(\Gamma)$$
se ramènera ainsi à celle de 
$$\frac{1}{d^j} \int \sum \int_{f^j(E_2(i,j) \cap Q_{\lambda} \cap \Gamma)} T \wedge
     [D_{\alpha}(\Gamma)] d \nu_{\Qcal}(\Gamma).$$
Ensuite les disques
     $D_{\alpha}(\Gamma)$ seront séparés en deux catégories: ceux qui
     sont proches de la direction stable et les autres. La majoration
     de 
$$\frac{1}{d^j}
\int \sum \int_{f^j(E_2(i,j) \cap Q_{\lambda} \cap \Gamma)} T \wedge
     [D_{\alpha}(\Gamma)] d \nu_{\Qcal}(\Gamma)$$
où la somme sera
     prise sur chacune des catégories sera alors effectuée dans le
     deuxième et le troisième paragraphe.

\subparagraph{{1) Approximation de $\int \frac{[f^j(\Gamma)]}{d^j} d \nu_{\Qcal}(\Gamma)$}}

$\mbox{ }$

Voici le plan de ce paragraphe: dans un premier
temps on va faire un découpage du but, c'est-à-dire on va quadriller
un ouvert $U$ qui contient essentiellement $f^j(\Gamma)$ en $4$-cubes
affines. Puis on complètera les courbes
$f^j(\Gamma)$ de sorte à enlever les problèmes liés au bord et enfin on construira un courant uniformément
géométrique qui approximera bien $ \int \frac{[f^j(\Gamma)]}{d^j} d
\nu_{\Qcal}(\Gamma)$.

\subparagraph{{i) Découpage du but}}

Quitte à considérer un recouvrement fini de $\Pp^2$ par des ouverts
$U$, on pourra supposer que $f^j(E_2(i,j)) \subset U$.

On considère maintenant un découpage $\Qcal^{'}$ d'un voisinage de $\overline{U}$
en $4$-cubes affines de taille $\frac{1}{l}$. Rappelons que
$\Qcal^{'}$ est de la forme
$$\Qcal^{'}= \{ \pi_1^{-1}(s_1^{'}) \times \pi_2^{-1}(s_2^{'}) \mbox{
  , } (s_1^{'},s_2^{'}) \in \Scal_1^{'} \times  \Scal_2^{'} \}$$
où $\Scal_1^{'}$ et $\Scal_2^{'}$ sont des quadrillages de $\Cc$ de
  taille $\frac{1}{l}$. Dans la suite ce quadrillage sera supposé bien
  choisi par rapport au courant $ \int \frac{[f^j(\Gamma)]}{d^j} d
\nu_{\Qcal}(\Gamma)$ et à la mesure $\nu$. Cela signifie que $T \wedge  \int \frac{[f^j(\Gamma)]}{d^j} d
\nu_{\Qcal}(\Gamma)$ et $\nu$ chargent peu un petit voisinage du bord
  de $\Qcal^{'}$ (voir le lemme 4.5 de \cite{Du2}). La position du
  quadrillage $\Qcal^{'}$ dépend donc de $j$ mais sa taille
  (i.e. $\frac{1}{l}$) sera toujours indépendante de $j$. Si on note
  $\Qcal^{'}_{\lambda^{'}}=\cup Q^{'}_{\lambda^{'}}$ où les
  $Q^{'}_{\lambda^{'}}$ sont les homothétiques des $Q^{'} \in
  \Qcal^{'}$ de rapport $\lambda^{'}$, on pourra supposer que
  $f^j(E_2(i,j)) \subset \Qcal^{'}_{\lambda^{'}}$ quitte à remplacer
  $E_2(i,j)$ par $E_2(i,j) \cap f^{-j}(\Qcal^{'}_{\lambda^{'}})$. Ici $\lambda^{'}$ ne dépend que de $\nu$.

\subparagraph{{ii) Complétion de $f^j(\Gamma)$}}

Il s'agit ici d'enlever le problème que crée le bord de $f^j(\Gamma)$
dans la construction des disques dans cette courbe.

Remarquons tout d'abord que l'on peut enlever les $\Gamma$ du support
de $\nu_{\Qcal}$ pour lesquels $\mbox{Aire}(f^j(\Gamma))$ est
supérieure à $Kd^j$. En effet, si $\Gcal$ désigne l'ensemble de ces
$\Gamma$, on a:
$$K d^j \nu_{\Qcal}(\Gcal) \leq \int \int [f^j(\Gamma)] \wedge \omega
d \nu_{\Qcal}(\Gamma) \leq \int f_{*}^{j} S_{\Qcal} \wedge \omega \leq
 \int f_{*}^{j} S \wedge \omega = d^j.$$
Autrement dit, $\nu_{\Qcal}(\Gcal) \leq \frac{1}{K}$ et le fait
d'enlever ces graphes $\Gamma$ à $T \wedge S_{\Qcal}$ modifie peu
cette mesure.

Ensuite, quitte à changer $Q$ en un cube un peu plus petit, nous
pouvons supposer que la longueur de $\partial f^j(\Gamma)$ est majorée
par $K d^{j/2}$ (en utilisant un argument longueur-aire).

$ $

Expliquons maintenant comment on complète les courbes $f^j(\Gamma)$.

Considérons un carré $s_1^{'}$ de $\Scal_1^{'}$ et un graphe
$\Delta^{'}$ d'une
fonction holomorphe au-dessus
de celui-ci d'aire inférieure à $\frac{1}{2}$. Supposons que l'intersection entre $\Delta^{'}$ et $f^j(\Gamma)$ contienne un ouvert
de $\Delta^{'}$. En utilisant une inégalité isopérimétrique, on constate que $\Delta^{'}$ vérifie un des cas suivants:

- Soit $\mbox{Longueur}(\pi_1( \partial f^j(\Gamma) \cap \Delta^{'}))
  \geq \frac{(1- \lambda^{'}) \delta^{'}}{l}$. Ici $\delta^{'}$ est
  une constante petite qui sera précisée plus loin.

- Soit cette longueur est plus petite que $\frac{(1- \lambda^{'})
  \delta^{'}}{l}$ mais l'aire de $\pi_1(f^j(\Gamma) \cap \Delta^{'})$ est
supérieure à $\frac{1}{l^2}-\frac{(1- \lambda^{'})^2
  \delta^{'2}}{l^2}$.

- Soit $f^j(\Gamma) \cap \Delta^{'}$ est confinée près du bord de
$\Qcal^{'}$ (i.e. $f^j(\Gamma) \cap \Delta^{'} \subset
(\Qcal^{'}_{\lambda^{''}})^c$ avec $1 > \lambda^{''} > \lambda^{'}$).

- Soit $f^j(\Gamma) \Subset \Delta^{'}$ n'est pas confinée près du
bord et alors la masse pour $T$ de
$f^j(\Gamma)$  est négligeable (si $\delta^{'}$ est petit). Quitte à
enlever ces $\Gamma$ à $T \wedge S_{\Qcal}$ (ce qui ne change presque
pas cette mesure), on pourra supposer
qu'aucun $\Gamma$ ne vérifie ce dernier cas.

Ce sont les composantes $\Delta^{'}$ des deux premières catégories que
l'on ajoute à $f^j(\Gamma)$.

En faisant la même chose avec les autres carrés $s_1^{'}$ de
$\Scal_1^{'}$ puis avec $\Scal_2^{'}$, on obtient ainsi un prolongement
$\widetilde{f^j(\Gamma)}$ de $f^j(\Gamma)$.

\begin{remarque*}

A ce stade de la démonstration $\widetilde{f^j(\Gamma)}$ peut être égal à
$f^j(\Gamma)$. En effet rien ne dit que l'on peut effectivement
compléter $f^j(\Gamma)$.

Dans la suite nous verrons que les graphes $\Delta^{'}$ que nous
arriverons à ajouter à $f^j(\Gamma)$ seront des limites de graphes contenus dans des courbes du type $f^p(L)$.

\end{remarque*}

\subparagraph{{iii) Construction du courant uniformément géométrique}}

La construction est la même que celle décrite après la proposition
\ref{prop4} du paragraphe \ref{geometrique}.

Dans les composantes connexes de $\widetilde{f^j(\Gamma)}$ au-dessus d'un carré $s_i^{'}$ de
$\Scal_i^{'}$ (i.e. dans $\pi_i^{-1}(s_i^{'}) \cap \widetilde{f^j(\Gamma)}$), on a un certain
nombre de disques qui sont des graphes au-dessus de $s_i^{'}$ (qui peuvent éventuellement s'intersecter entre eux). En sommant les
courants d'intégration sur les graphes qui ont une aire inférieure à
$\frac{1}{2}$, on obtient ainsi deux courants inférieurs à
$\widetilde{f^j(\Gamma)}$: $[C_{\Qcal_1^{'}}^{'}]$ et
$[C_{\Qcal_2^{'}}^{'}]$ où les $C_{\Qcal_i^{'}}^{'}$ sont des courbes
à bord dans $\pi_i^{-1}(\partial \Scal_i^{'})$. En réunissant
maintenant ces deux courbes, on obtient:
$$[C_{\Qcal_1^{'}}^{'} \cup C_{\Qcal_2^{'}}^{'}]=\sum
  [D_{\alpha}(\Gamma)],$$
où les $D_{\alpha}(\Gamma)$ sont des disques dans les cubes $Q^{'}$ de
  $\Qcal^{'}$. Le courant $R_j= \frac{1}{d^j}
\int \sum [D_{\alpha}(\Gamma)] d \nu_{\Qcal}(\Gamma)$ ainsi construit
est uniformément géométrique dans les cubes $Q^{'}$ de $\Qcal^{'}$ et
il vérifie:

\begin{lemme}

$$\langle \frac{1}{d^j} \int [\widetilde{f^j(\Gamma)} \cap
  \Qcal_{\lambda^{'}}^{'}] d \nu_{\Qcal}(\Gamma) -
  R_{j|\Qcal_{\lambda^{'}}^{'}}, \pi_1^{*} \omega + \pi_2^{*} \omega
  \rangle \leq \frac{C}{l^2}$$
Ici $C$ est indépendante de $j$.
\end{lemme}

\begin{proof}

Considérons les courants $S_{m, \Qcal}= \int [\Gamma] d \nu_{m,
  \Qcal}(\Gamma)$ qui approximent le courant $S$ (voir le paragraphe
  \ref{geometrique}).

L'inégalité ci-dessus est alors vraie en remplaçant $\nu_{\Qcal}$ par
  $\nu_{m, \Qcal}$ et $R_j$ par $R_{j,m}$ (où $R_{j,m}$ est construit
  de façon analogue à $R_j$ à partir des courants $\frac{f_{*}^{j}
  S_{m, \Qcal}}{d^j}$). Cela provient, exactement comme dans la
  proposition \ref{prop4}, de la formule de Riemann-Hurwitz
  et du fait que le genre d'une courbe $f^p(L)$ ($p
  \in \Nn$) vaut $0$. Par ailleurs dans cette nouvelle inégalité le
  $C$ est indépendant de $m$ et de $j$.

Maintenant, l'application $F(\Gamma)=\langle [\widetilde{f^j(\Gamma)} \cap
  \Qcal_{\lambda^{'}}^{'}] - \sum [D_{\alpha}(\Gamma) \cap \Qcal_{\lambda^{'}}^{'}],\pi_1^{*} \omega + \pi_2^{*} \omega
  \rangle$ est semi-continue inférieurement.

On obtient donc le lemme en passant à la limite sur les mesures
  $\nu_{m, \Qcal}$.

\end{proof}

Maintenant, si on applique $T$, on obtient que la différence
$$\frac{1}{d^j} \int \int_{\Qcal_{\lambda^{'}}^{'}} T \wedge
[\widetilde{f^j(\Gamma)}]  d \nu_{\Qcal}(\Gamma) -
\int_{\Qcal_{\lambda^{'}}^{'}} T \wedge R_j$$
est majorée par $\epsilon(l)$ qui est une suite qui tend vers $0$ quand
$l$ croît vers l'infini (voir le paragraphe \ref{courant}). Cela
résulte d'une part du lemme précédent et d'autre part qu'étant donné que le
quadrillage $\Qcal^{'}$ est supposé bien choisi par rapport à $T \wedge  \int \frac{[f^j(\Gamma)]}{d^j} d
\nu_{\Qcal}(\Gamma)$ il l'est pour
$$\frac{1}{d^j} \int  T \wedge [\widetilde{f^j(\Gamma)}]  d \nu_{\Qcal}(\Gamma) -
\int T \wedge R_j$$
(voir le lemme 4.5 de \cite{Du2}).

Par ailleurs, comme nous l'avons dit dans le paragraphe
\ref{geometrique}, la suite $\epsilon(l)$ est indépendante de $j$.

Pour majorer la mesure de $E_2$, nous sommes donc ramenés à majorer 
$$\frac{1}{d^j} \int \int_{f^j(E_2(i,j) \cap \Gamma \cap Q_{\lambda})} \sum
T \wedge [D_{\alpha}(\Gamma)]  d \nu_{\Qcal}(\Gamma)$$
par $\epsilon$.

$ $

Ce sont maintenant les disques $D_{\alpha}(\Gamma)$ que nous allons
séparer en deux catégories. Par ailleurs, comme $f^j(E_2(i,j)) \subset
\Qcal^{'}_{\lambda^{'}}$, on ne considèrera dans la suite que les disques $D_{\alpha}(\Gamma)$ qui
entrent dans $\Qcal^{'}_{\lambda^{'}}$.

Fixons $\delta > 0$ et $Q^{'}$ un cube de $\Qcal^{'}$. Nous pouvons
supposer que la variation de la pente de $D_{\alpha}(\Gamma)$ dans
$Q_{\lambda^{'}}^{'}$ est inférieure à $\frac{\delta}{4}$ (quitte à
redécouper les disques $D_{\alpha}(\Gamma)$ via un sous-quadrillage de
$\Qcal^{'})$. De plus, la fonction $x \mapsto E_S(x)$ étant continue
sur $f^j(E_2(i,j))$, nous pouvons aussi supposer que sur les
$Q_{\lambda^{'}}^{'}$, la variation de cette fonction est inférieure à
$\frac{\delta}{4}$.

Maintenant, si $D_{\alpha}(\Gamma)$ est un disque de
$\widetilde{f^j(\Gamma)}$ tel que $\int_{f^j(E_2(i,j) \cap \Gamma
  \cap Q_{\lambda})} T \wedge [D_{\alpha}(\Gamma)] > 0$ alors
nécessairement $D_{\alpha}(\Gamma)$ contient des points de
$\pi(\widehat{E}_2)$. La séparation des disques $D_{\alpha}(\Gamma)$ se fait de la façon
suivante. Nous noterons $\Dcal_1$ l'ensemble des disques qui
contiennent un point $x$ de $\pi(\widehat{E}_2)$ avec l'angle entre la tangente à $D_{\alpha}(\Gamma)$ en $x$ et
$E_S(x)$ supérieur à $\delta$. Les autres disques $D_{\alpha}(\Gamma)$
(qui sont donc proches de la direction stable) constituent l'ensemble
$\Dcal_2$. Dans la suite, nous aurons donc deux majorations à
effectuer: celle de 
$$\frac{1}{d^j} \int \int_{f^j(E_2(i,j) \cap \Gamma \cap Q_{\lambda})} \sum_{\Dcal_1}
T \wedge [D_{\alpha}(\Gamma)]  d \nu_{\Qcal}(\Gamma)$$
et celle de 
$$\frac{1}{d^j} \int \int_{f^j(E_2(i,j) \cap \Gamma \cap Q_{\lambda})} \sum_{\Dcal_2}
T \wedge [D_{\alpha}(\Gamma)]  d \nu_{\Qcal}(\Gamma)$$
par $\epsilon$.

Ces majorations vont être effectuées dans les deux paragraphes suivants.

\subparagraph{{2) Majoration de $\frac{1}{d^j} \int \int_{f^j(E_2(i,j) \cap \Gamma \cap Q_{\lambda})} \sum_{\Dcal_1}
T \wedge [D_{\alpha}(\Gamma)]  d \nu_{\Qcal}(\Gamma)$}}

$\mbox{ }$

Nous savons que le nombre de disques $D_{\alpha}(\Gamma)$ dans $\widetilde{f^j(\Gamma)}$
est majoré par $C(l)d^j$. En effet les disques que l'on a ajoutés à
$f^j(\Gamma)$ pour obtenir $\widetilde{f^j(\Gamma)}$ contiennent soit
un morceau consistant du bord de $f^j(\Gamma)$ (i.e. de longueur
supérieure à  $\frac{(1- \lambda^{'}) \delta^{'}}{l}$), soit une aire
minorée par $\frac{1}{l^2}-\frac{(1- \lambda^{'})^2
  \delta^{'2}}{l^2}$. Afin de
simplifier les expressions, nous supposerons dans la suite que $C(l) \leq 1$. De
la même façon, la mesure $\nu_{\Qcal}$ sera considérée de masse $1$.

Nous allons maintenant faire la majoration de $\frac{1}{d^j} \int \int_{f^j(E_2(i,j) \cap \Gamma \cap Q_{\lambda})} \sum_{\Dcal_1}
T \wedge [D_{\alpha}(\Gamma)]  d \nu_{\Qcal}(\Gamma)$, en admettant
certaint points que nous démontrerons ensuite.

On découpe $\Pp^2$ en ouverts disjoints simplement connexes $B$ tels que $\tau_x^{-1}(B) \subset
B(0, \frac{1}{4 p^4})$ pour tout $x$ de $B$. Ensuite, nous noterons
$\cup \mbox{BC}$ l'union des bonnes composantes dans les $f^{-i}(B)$
(i.e. des ouverts de la forme $\tau_{y_{-i}} \circ f_{\widehat{y}}^{-i} \circ \tau_{y_0}^{-1} (B)$ avec $\widehat{y}$ et
$\widehat{f}^{-i}(\widehat{y})$ dans $\widehat{E}_2$ et $y_0$ dans
$B$).

Remarquons que par définition de $E_2(i,j)$, nous pouvons recouvrir $f^j(E_2(i,j))$ par des bonnes
composantes $\tau_{y_{-i}} \circ f_{\widehat{y}}^{-i} \circ
\tau_{y_0}^{-1} (B)$ avec $y_{-i} \in f^j(E_2(i,j))$. Par ailleurs, si $i$ est assez grand, celles-ci évitent
un petit voisinage du support de $\mu$. En effet, d'une part
$d(\pi(\widehat{x}), \mbox{Support}(\mu)) \geq \frac{1}{p}$ si
$\widehat{x}$ est dans $\widehat{E}_2$ (voir le paragraphe \ref{a}) et d'autre
part, les bonnes composantes sont de diamètre exponentiellement petit
(voir le point $2$ du théorème du paragraphe \ref{inverse}). Le nombre
des bonnes composantes qui recouvrent $f^j(E_2(i,j))$ est donc majoré
par $Nd^i$ (voir le paragraphe $\gamma$ de la section \ref{nbe}).

Nous verrons que si $\mbox{BC}$ est une de ces bonnes composantes
alors $\int_{\mbox{BC}} T \wedge [D_{\alpha}(\Gamma)] \leq C
d^{-i}$. Cela viendra du fait que $D_{\alpha}(\Gamma)$ n'est pas
proche de la direction stable. Maintenant, si $\frac{1}{d^j} \int \int_{f^j(E_2(i,j) \cap \Gamma \cap Q_{\lambda})} \sum_{\Dcal_1}
 T \wedge [D_{\alpha}(\Gamma)]  d \nu_{\Qcal}(\Gamma)$ est supérieur à
$\epsilon$, alors il existe un disque $D_{\alpha}(\Gamma)$ avec
$\int_{f^j(E_2(i,j) \cap \Gamma \cap Q_{\lambda})} T \wedge
[D_{\alpha}(\Gamma)] \geq \epsilon$. Par ce qui précède, ce disque coupe
au moins $\frac{\epsilon}{C} d^i$ bonnes composantes (car celles-ci
sont disjointes). De plus, ces
bonnes composantes sont incluses dans un $(e^{- \beta i + \gamma
  i})$-voisinage $V$ de ce disque (par le point $2$ du théorème de la
section \ref{inverse}). Nous verrons que la mesure $\epsilon(i)$ pour
$\nu$ de ce voisinage tend vers $0$ quand $i$ croît vers l'infini et
ceci uniformément sur les disques $D_{\alpha}(\Gamma)$.

Quitte à remplacer $E_2(i,j)$ par $E_2(i,j) \setminus f^{-j}(V)$
nous sommes ramenés à majorer 
$$\frac{1}{d^j} \int \int_{f^j(E_2(i,j) \cap \Gamma \cap Q_{\lambda})} \sum_{\Dcal_1}
T \wedge [D_{\alpha}(\Gamma)]  d \nu_{\Qcal}(\Gamma)$$
par $\epsilon - \epsilon(i)$ avec $f^j(E_2(i,j))$ qui évite $V$.

Comme précédemment, si ce n'est pas la cas, on obtient un disque
$D_{\alpha}(\Gamma)$ avec 
$$\int_{f^j(E_2(i,j) \cap \Gamma \cap Q_{\lambda})} T \wedge
[D_{\alpha}(\Gamma)] \geq \epsilon - \epsilon(i).$$
Ce disque coupe donc
au moins $\frac{\epsilon - \epsilon(i)}{C}d ^i$ bonnes composantes qui
sont disjointes des précédentes.

On a donc construit $\left( \frac{\epsilon}{C}+\frac{\epsilon -
  \epsilon(i)}{C} \right) d ^i$ bonnes composantes hors d'un petit
voisinage du support de $\mu$.

En itérant maintenant le procédé et en utilisant le fait qu'il y a au
plus $Nd^i$ bonnes composantes hors de ce voisinage du support de
$\mu$, on obtient l'existence d'un $k \leq \frac{2CN}{\epsilon}$ pour
lequel $\int_{f^j(E_2(i,j) \cap \Gamma \cap Q_{\lambda})} T \wedge
[D_{\alpha}(\Gamma)] \leq \epsilon -k \epsilon(i)$ (avec $i$ grand). Ici
$E_2(i,j)$ est le $E_2(i,j)$ initial auquel on a enlevé $k$
voisinages du même type que $V$.

Nous avons donc la majoration que nous voulions.

$ $

Il reste maintenant à montrer les points que nous avons admis: le fait
que $\int_{\mbox{BC}} T \wedge [D_{\alpha}(\Gamma)] \leq C
d^{-i}$ où $\mbox{BC}$ est une bonne composante de $f^{-i}(B)$ et la
majoration de la mesure pour $\nu$ d'un $(e^{- \beta i + \gamma
  i})$-voisinage de $D_{\alpha}(\Gamma)$ par $\epsilon(i)$.

\subparagraph{{i) Majoration de $\int_{\mbox{BC}} T \wedge
    [D_{\alpha}(\Gamma)]$ par $ Cd^{-i}$}}

Rappelons d'une part que l'on considère ici une bonne composante $\mbox{BC}=\tau_{y_{-i}} \circ
f_{\widehat{y}}^{-i} \circ \tau_{y_0}^{-1} (B)$ (avec $\widehat{y}$ et
$\widehat{f}^{-i}(\widehat{y})$ dans $\widehat{E}_2$, $y_0$ dans $B$
et $y_{-i} \in f^j(E_2(i,j))$), et d'autre part que $\tau_{y_0}^{-1}(B)
\subset B(0,\frac{1}{4p^4})$ et $\frac{1}{p} \leq \|
C_{\gamma}^{\pm 1}(\widehat{f}^{-i}(\widehat{y})) \| \leq p$.

Comme dans la majoration de la mesure de $E_1$, si on note
$D_{\alpha}^0(\Gamma)$ la partie de
$C_{\gamma}^{-1}(\widehat{f}^{-i}(\widehat{y})) \circ
\tau_{y_{-i}}^{-1}(D_{\alpha}(\Gamma))$ susceptible de rencontrer
$C_{\gamma}^{-1}(\widehat{f}^{-i}(\widehat{y})) \circ
f_{\widehat{y}}^{-i} B(0, \frac{1}{p})$, alors $D_{\alpha}^0(\Gamma)$
se met sous la forme d'un graphe $(z, \psi(z))$ où $z$ est dans un
disque $D_0$ de taille $e^{-\chi_1(\widehat{y})i + \gamma i}$ et $|
\psi^{'}(z) | \leq e^{-\chi_1(\widehat{y})i + \gamma i}$ sur $D_0$
(on s'est placé dans un repère orthonormé où l'axe des abscisses est parallèle à une tangente à $D_{\alpha}^0(\Gamma)$ en
un point $x$ de $C_{\gamma}^{-1}(\widehat{f}^{-i}(\widehat{y})) \circ
f_{\widehat{y}}^{-i} \circ \tau_{y_0}^{-1} (B)$).

Maintenant, comme l'angle entre $D_{\alpha}(\Gamma)$ et
$E_S(\widehat{f}^{-i}(\widehat{y}))=D_0 \tau_{y_{-i}}
E_1(\widehat{f}^{-i}(\widehat{y}))$ est supérieur à
$\frac{\delta}{2}$, la partie $\Delta_{\alpha}^0(\Gamma)$ de
$D_{\alpha}^0(\Gamma)$ susceptible de rencontrer $C_{\gamma}^{-1}(\widehat{f}^{-i}(\widehat{y})) \circ
f_{\widehat{y}}^{-i} B(0, \frac{1}{p})$ se met aussi sous la forme
d'un graphe $(z, \phi(z))$ dans le repère orthonormé
$$(C_{\gamma}^{-1}(\widehat{f}^{-i}(\widehat{y}))E_2(\widehat{f}^{-i}(\widehat{y})),C_{\gamma}^{-1}(\widehat{f}^{-i}(\widehat{y}))E_1(\widehat{f}^{-i}(\widehat{y}))).$$
Ici $z$ est dans un disque $D$ de taille essentiellement
$e^{-\chi_1(\widehat{y})i + \gamma i}$ et $| \phi^{'}(z) | \leq
C(\delta)$ pour $z \in D$.

Nous pouvons maintenant utiliser la transformée de graphe. Plus
précisément, si on note comme dans le paragraphe \ref{exposant},
$g_{\widehat{x}}=C_{\gamma}^{-1}(\widehat{f}(\widehat{x})) \circ f_x
\circ C_{\gamma}(\widehat{x})$ (où $\pi(\widehat{x})=x$) on a la

\begin{proposition}

Soit $\gamma_0 > 0$. Si $i$ est assez grand, l'image de $\Delta_{\alpha}^0(\Gamma) \cap
C_{\gamma}^{-1}(\widehat{f}^{-i}(\widehat{y})) \circ f_{\widehat{y}}^{-i} B(0,
\frac{1}{p^4})$ par $g_{\widehat{f}^{-1}(\widehat{y})} \circ \cdots
  \circ g_{\widehat{f}^{-i}(\widehat{y})}=C_{\gamma}^{-1}(\widehat{y})
    \circ \tau_{y_0}^{-1} \circ f^i \circ \tau_{y_{-i}} \circ
      C_{\gamma}(\widehat{f}^{-i}(\widehat{y}))$ est un graphe
      $(z,\psi(z))$ au-dessus d'une partie $D_i$ de $C_{\gamma}^{-1}(\widehat{y})
      E_2(\widehat{y})$ qui vérifie $| \psi(z_1) - \psi(z_2) | \leq
      \gamma_0 |z_1 - z_2|$ si $z_1$ et $z_2$ sont dans $D_i$.

\end{proposition}

\begin{proof}

La démonstration va utiliser la généralisation de la transformée de
graphe qui se trouve dans l'appendice.

Les fonctions $g_{\widehat{f}^{-l}(\widehat{y})}$ ont les propriétés
suivantes (voir le paragraphe \ref{exposant}):

- $D_0 g_{\widehat{f}^{-l}(\widehat{y})}=\left( \begin{array}{cc}
\mu_1(\widehat{f}^{-l}(\widehat{y})) & 0\\
0 & \mu_2(\widehat{f}^{-l}(\widehat{y}))\\
\end{array} \right) $
avec pour $i=1,2$
$$e^{\chi_i (\widehat{y})- \gamma} \leq |\mu_i(\widehat{f}^{-l}(\widehat{y}))|
\leq e^{\chi_i (\widehat{y})+ \gamma}.$$
- Sur $B(0,\delta(\widehat{f}^{-l}(\widehat{y})))$, on a
$$g_{\widehat{f}^{-l}(\widehat{y})}(w)=D_0
g_{\widehat{f}^{-l}(\widehat{y})}w+h(w)$$
avec $\| D_w h \| \leq \frac{\gamma \| w
  \|}{\delta(\widehat{f}^{-l}(\widehat{y}))}$.

Autrement dit, dans le repère 
$$(C_{\gamma}^{-1}(\widehat{f}^{-l}(\widehat{y}))E_1(\widehat{f}^{-l}(\widehat{y})),C_{\gamma}^{-1}(\widehat{f}^{-l}(\widehat{y}))E_2(\widehat{f}^{-l}(\widehat{y}))),$$
$g_{\widehat{f}^{-l}(\widehat{y})}$ s'écrit:
$$g_{\widehat{f}^{-l}(\widehat{y})}=(\mu_1(\widehat{f}^{-l}(\widehat{y}))x+
\alpha(x,y),\mu_2(\widehat{f}^{-l}(\widehat{y}))y + \beta(x,y))$$
où $| \alpha |_{C^1}, | \beta |_{C^1} \leq (e^{- \chi_1
  (\widehat{y})l+3 \gamma l}) \frac{\gamma}{p}$ sur $B(0,\frac{e^{- \chi_1
  (\widehat{y})l+2 \gamma l}}{p^2})$ (en effet
$\delta(\widehat{f}^{-l}(\widehat{y})) \geq \frac{e^{- \gamma l}}{p}$
  puisque $\delta$ est une fonction tempérée et $\delta(\widehat{y})
  \geq \frac{1}{p}$).

Nous considérons $| \alpha |_{C^1}$ et $| \beta |_{C^1}$ sur $B(0,\frac{e^{- \chi_1
  (\widehat{y})l+2 \gamma l}}{p^2})$ car cette boule contient $C_{\gamma}^{-1}(\widehat{f}^{-l}(\widehat{y})) \circ
f_{\widehat{y}}^{-l} B(0,\frac{1}{p^4})$ (voir le point $2$ du
théorème section \ref{inverse}).

Commençons maintenant par regarder l'image de $\Delta_{\alpha}^0(\Gamma) \cap
C_{\gamma}^{-1}(\widehat{f}^{-i}(\widehat{y})) \circ f_{\widehat{y}}^{-i} B(0,
\frac{1}{p^4})$ par $g_{\widehat{f}^{-i}(\widehat{y})}$. Par
l'appendice cette image est un graphe au-dessus d'un certain ouvert de
$$C_{\gamma}^{-1}(\widehat{f}^{-i+1}(\widehat{y}))E_2(\widehat{f}^{-i+1}(\widehat{y})).$$
Il s'écrit donc, dans le repère
$$(C_{\gamma}^{-1}(\widehat{f}^{-i+1}(\widehat{y}))E_1(\widehat{f}^{-i+1}(\widehat{y})),C_{\gamma}^{-1}(\widehat{f}^{-i+1}(\widehat{y}))E_2(\widehat{f}^{-i+1}(\widehat{y}))),$$
$(\phi_1(z),z)$ avec $z \in D_1$. Par ailleurs,
toujours grâce à l'appendice, on a pour $z_1$ et $z_2$ dans $D_1$,
$$| \phi_1(z_1) - \phi_1(z_2) | \leq
\frac{\mu_1(\widehat{f}^{-i}(\widehat{y}))C(\delta) + \max(| \alpha |_{C^1},| \beta |_{C^1})
  (1+C(\delta))}{\mu_2(\widehat{f}^{-i}(\widehat{y}))- \max(| \alpha |_{C^1},| \beta |_{C^1})
  (1+C(\delta))} | z_1 - z_2 |$$
c'est-à-dire
$$| \phi_1(z_1) - \phi_1(z_2) | \leq \left(
\frac{\mu_1(\widehat{f}^{-i}(\widehat{y}))}{\mu_2(\widehat{f}^{-i}(\widehat{y}))}
+ \gamma \right) C(\delta) | z_1 - z_2 | \leq t C(\delta) | z_1 - z_2
|$$
avec $0 < t < 1$ le tout si $i$ est grand.

Enfin, on peut remarquer que le graphe que l'on vient de construire
est dans $C_{\gamma}^{-1}(\widehat{f}^{-i+1}(\widehat{y})) \circ
f_{\widehat{y}}^{-i+1} B(0,\frac{1}{p^4})$.

Si $i$ est grand, on peut donc recommencer ce que l'on vient de faire
$i_0$ fois et on obtient ainsi un graphe $(\phi_{i_0}(z),z)$
au-dessus de $D_{i_0} \subset
C_{\gamma}^{-1}(\widehat{f}^{-i+i_0}(\widehat{y}))E_2(\widehat{f}^{-i+i_0}(\widehat{y}))$
qui vérifie:
$$| \phi_{i_0}(z_1) - \phi_{i_0}(z_2) | \leq t^{i_0}C(\delta) | z_1 - z_2
| \leq \gamma_0 | z_1 - z_2|,$$
si $z_1$ et $z_2$ sont dans $D_{i_0}$.

Ensuite, toujours grâce à la transformée de graphe (voir l'appendice)
et à la majoration des termes d'ordre deux de
$g_{\widehat{f}^{-l}(\widehat{y})}$ (i.e. $| \alpha |_{C^1}$ et $| \beta |_{C^1}$) par $(e^{- \chi_1
  (\widehat{y})l+3 \gamma l}) \frac{\gamma}{p}$ sur $B(0,\frac{e^{- \chi_1
  (\widehat{y})l+2 \gamma l}}{p^2})$, quand on continue à pousser en
avant le graphe $(\phi_{i_0}(z),z)$, on obtient
toujours des graphes $\phi_l$ au-dessus d'une région $D_l$ de
$C_{\gamma}^{-1}(\widehat{f}^{-i+l}(\widehat{y}))E_2(\widehat{f}^{-i+l}(\widehat{y}))$.
De plus, ces graphes vérifient:
$$| \phi_{l}(z_1) - \phi_{l}(z_2) | \leq \gamma_0 | z_1 - z_2|,$$
pour $z_1$ et $z_2$ dans $D_l$ et $l$ compris entre $i_0 +1$ et $i$.

Le cran $i$ donne alors ce que l'on voulait démontrer.

\end{proof}

Notons $g$ le graphe obtenu par la proposition.

Il nous reste à calculer 
$$\int_{\mbox{BC}} T \wedge [D_{\alpha}(\Gamma)]= \frac{1}{d^i} \int_B T
\wedge  [ \tau_{y_0} \circ C_{\gamma}(\widehat{y})(g)] \leq \frac{1}{d^i}
\int_{C_{\gamma}^{-1}(\widehat{y})B(0,\frac{1}{4p^4})}
C_{\gamma}^{*}(\widehat{y}) \tau_{y_0}^{*} T \wedge [g].$$
Comme dans le paragraphe précédent, nous allons utiliser un argument
de capacité.

En effet, si on note
$\Bcal_l=C_{\gamma}^{-1}(\widehat{y})B(0,\frac{1}{lp^4})$
($l=1,2,3,4$), $\Delta$ une direction complexe parallèle à
$C_{\gamma}^{-1}(\widehat{y})E_2(\widehat{y})$ qui coupe $g$ en un
point de $\Bcal_4$ et $\pi$ la projection
orthogonale sur $\Delta$, on a d'une part $\pi(g \cap \Bcal_4) \subset
\Delta \cap \Bcal_3$ et d'autre part $\Delta \cap \Bcal_2 \subset
\pi(g \cap \Bcal_1)$ (si $\gamma_0$ est assez petit par rapport à une
constante qui dépend de $p$).

Autrement dit, comme dans le paragraphe précédent,
$$\int_{\mbox{BC}} T \wedge [D_{\alpha}(\Gamma)] \leq \frac{C}{d^i}
\mbox{Cap}(\Delta \cap \Bcal_3, \Delta \cap \Bcal_2) \leq
\frac{C}{d^i}.$$

Il nous reste maintenant à majorer la mesure pour $\nu$ d'un $(e^{-
  \beta i + \gamma i})$-voisinage d'un disque $D_{\alpha}(\Gamma)$ par
  une suite $\epsilon(i)$ qui tend vers $0$ uniformément sur les $D_{\alpha}(\Gamma)$.

\subparagraph{{ii) Majoration de $\nu((e^{-
  \beta i + \gamma i} ) \mbox{-voisinage d'un disque }
  D_{\alpha}(\Gamma))$}}

Il s'agit ici de voir que
$$ \forall \alpha \mbox{ , } \exists i_0 \mbox{ , } \forall i \geq i_0
\mbox{ , } \forall D_{\alpha}(\Gamma) \mbox{ , } \nu((e^{-
  \beta i + \gamma i} ) \mbox{-voisinage d'un disque }
  D_{\alpha}(\Gamma)) \leq \alpha.$$
Faisons un raisonnement par l'absurde. Si ce qui précède est faux, on
obtient une suite de disques $D_n$ (d'aire inférieure à $\frac{1}{2}$)
et une suite $i_n$ telles que $ \nu((e^{-
  \beta i_n + \gamma i_n} ) \mbox{-voisinage de }
  D_n) \geq \alpha$.

Quitte à extraire une sous-suite, $D_n$ converge vers un disque $D$ et
on a nécessairement $\nu(D) \geq \frac{\alpha}{2}$.

Pour obtenir une contradiction, nous allons montrer que $D$ est inclus
dans $\{ x \mbox{, } \nu(S,x) > 0 \}$, ce qui contredira le fait que
$\nu$ ne charge pas les ensembles algébriques.

Le disque $D$ vit dans un cube $Q^{'}$ de $\Qcal^{'}$. On peut
construire un quadrillage $\Qcal^{''}$ d'un voisinage $V$ de
$\overline{Q^{'}}$ en $4$-cubes affines et un courant $S_{\Qcal^{''}}$
uniformément géométrique dans ces cubes tels que la différence $(T
\wedge S - T \wedge S_{\Qcal^{''}})(V)$ soit très petite devant
$\alpha$.

Quitte à supposer que $\nu$ ne charge pas $\partial Q^{'}$ ainsi que
le bord de $\Qcal^{''}$, on peut trouver $Q^{''} \in \Qcal^{''}$ avec
$Q^{''} \Subset Q^{'}$ tel que:
$$\int \int_{D \cap Q^{''}_{\lambda^{''}}} T \wedge [\Gamma] d
  \nu_{\Qcal^{''}}(\Gamma) \geq \epsilon(\alpha) > 0$$
où $ Q^{''}_{\lambda^{''}}$ est l'homothétique de $Q^{''}$ de rapport
  $\lambda^{''}$ et $S_{\Qcal^{''}}= \int [\Gamma] d
  \nu_{\Qcal^{''}}(\Gamma)$.

Maintenant, si $\Gamma$ est un disque de $S_{\Qcal^{''}}$ dans
  $Q^{''}$, on a deux possibilités: soit $D \cap \Gamma \cap
  Q^{''}_{\lambda^{''}}$ contient un ouvert de $D$, soit $ \int_{D
  \cap Q^{''}_{\lambda^{''}}} T \wedge [ \Gamma] =0$ (car si
  l'intersection entre $D$ et $\Gamma$ dans $Q^{''}_{\lambda^{''}}$
  contient une infinité de points, on est dans le premier cas).

Autrement dit, comme il y a un nombre fini de $\Gamma$ pour lesquels
  $D \cap \Gamma \cap  Q^{''}_{\lambda^{''}}$ contient un ouvert de
  $D$ (car l'aire de $D$ est majorée par $\frac{1}{2}$), il en existe un avec $\nu_{\Qcal^{''}}(\Gamma) \geq \epsilon^{'}(\alpha)$. Mais alors
  $S_{\Qcal^{''}}$ (et donc $S$) a des nombres de Lelong strictement
  positifs sur $D \cap \Gamma$ (par le théorème de Lelong) donc $D
  \subset \{x \mbox{, } \nu(S,x) > 0 \}$ et on obtient ainsi une contradiction.

\subparagraph{{3) Majoration de $\frac{1}{d^j} \int \int_{f^j(E_2(i,j) \cap \Gamma \cap Q_{\lambda})} \sum_{\Dcal_2}
T \wedge [D_{\alpha}(\Gamma)]  d \nu_{\Qcal}(\Gamma)$}}

$\mbox{ }$

Pour simplifier, nous supposerons toujours que le nombre de disques
$D_{\alpha}(\Gamma)$ dans $f^j(\Gamma)$ est majoré par $d^j$ et que
$\nu_{\Qcal}$ est de masse $1$. Rappelons que $f^j(E_2(i,j)) \subset
\Qcal^{'}_{\lambda^{'}}$.

Si la quantité ci-dessus est minorée par $\epsilon$, on peut trouver
un cube $Q^{'}$ de $\Qcal^{'}$ et un $\Gamma$ pour lesquels on a:
$$\frac{1}{d^j} \int \int_{f^j(E_2(i,j) \cap \Gamma \cap
        Q_{\lambda}) \cap Q_{\lambda^{'}}^{'}} \sum_{\Dcal_2}
T \wedge [D_{\alpha}(\Gamma)]  \geq \frac{\epsilon}{C(l)}$$
où $C(l)$ est le nombre de cubes de $\Qcal^{'}$.

Maintenant, on peut construire un cube $C$ et un polydisque $P$ qui
        ont les propriétés suivantes:

- $C \subset P \Subset Q_{\lambda^{''}}^{'}$ (avec $1 > \lambda^{''} >
        \lambda^{'}$).

- $C$ et $P$ ont leur base parallèle à $E_s(x)$ où $x$ est un point de
        $f^j(E_2(i,j) \cap \Gamma \cap Q_{\lambda}) \cap C$.

- L'ensemble $\Dcal$ des disques $D_{\alpha}(\Gamma)$ dans $\Dcal_2$ qui passent par
        $C$ est constitué de graphes au-dessus de  $E_s(x)$ dans
        $P$.

- $\frac{1}{d^j} \int \sum_{\Dcal} \int_{f^j(E_2(i,j) \cap \Gamma \cap
        Q_{\lambda}) \cap C} T \wedge [D_{\alpha}(\Gamma)]  \geq
        \frac{\epsilon}{C(l)C(\delta, \lambda^{'})}$ (dans la suite
        nous remplacerons $\frac{\epsilon}{C(l)C(\delta,
        \lambda^{'})}$ par $\epsilon$).

Cela implique qu'il y a au moins $\frac{\epsilon}{C(l)}d^{j}$ disques
        $D_{\alpha}(\Gamma)$ tels que $\int_{f^j(E_2(i,j) \cap \Gamma \cap
        Q_{\lambda}) \cap C} T \wedge [D_{\alpha}(\Gamma)]  > 0$ (car
        $\int_{D_{\alpha}(\Gamma) \cap Q^{'}_{\lambda^{'}}} T \leq
        C(l)$). On notera toujours $\Dcal$ l'ensemble de ces disques.

Si $D$ est l'intersection avec $P$ d'un de ces disques de $\Dcal$, il
        contient une infinité de points de $f^j(E_2(i,j) \cap \Gamma \cap
        Q_{\lambda}) \cap C$. On peut donc construire une branche inverse
        $\tau_{x_0} \circ f_{\widehat{f}^{j}(\widehat{x})}^{-j} \circ
        \tau_{x_j}^{-1}$ qui envoie $D$ dans $\Gamma$ (ici
        $\widehat{f}^j(\widehat{x}) \in \widehat{E}_2$, $\pi \widehat{f}^j(\widehat{x}) = x_j
        \in f^j(E_2(i,j) \cap \Gamma \cap
        Q_{\lambda}) \cap D \cap C$ et $x_0 \in E_2(i,j) \cap \Gamma \cap
        Q_{\lambda}$).

Maintenant, si on utilise la transformée de graphe en temps négatif,
        on constate que $\tau_{x_0} \circ
        f_{\widehat{f}^j(\widehat{x})}^{-j} \circ \tau_{x_j}^{-1}
        (D_1)$ (où $D_1$ est un disque $D_{\alpha}(\Gamma)$ de $\Dcal$
        intersecté avec $P$) est un graphe au-dessus du disque $D(x_0,
        e^{- \chi_1(\widehat{x})j - \gamma j})$ dans $E_S(x_0)=D_0
        \tau_{x_0}E_1(\widehat{x})$ et que par ailleurs ce graphe est
        exponentiellement plat (si $\delta$ et $P$ sont suffisamment
        petits).

Si $D_1$ est différent de $D$, on en déduit donc que la préimage de
        $D_1$ par $\tau_{x_0} \circ
        f_{\widehat{f}^j(\widehat{x})}^{-j} \circ \tau_{x_j}^{-1}$ ne
        peut pas être dans $\Gamma$ (sinon cela contredirait le fait que
        $\Gamma$ est un graphe pour $j$ grand).

Autrement dit, il existe au moins $\frac{\epsilon}{C(l)}d^j$
        bonnes composantes distinctes de la forme $\tau_{x_0} \circ
        f_{\widehat{f}^j(\widehat{x})}^{-j} \circ
        \tau_{x_j}^{-1}(Q_{\lambda^{'}}^{'})$ qui coupent $\Gamma$.

Maintenant, comme au paragraphe précédent, on enlève un $(e^{- \beta j
        + \gamma j})$-voisinage $V$ de $\Gamma$ à $E_2(i,j)$ (ce
        voisinage est de masse $\epsilon(j)$ pour $\nu$) et on
        recommence ce que l'on vient de faire en remplaçant
        $E_2(i,j)$ par $E_2(i,j) \setminus V$, puis on itère le
        procédé.

Comme le nombre de bonnes composantes dans les préimages
        $f^{-j}(Q_{\lambda^{'}}^{'})$ est borné par $Nd^j$, le procédé
        finit par s'arrêter et on obtient la majoration voulue.

\section{\bf{Un exemple}}{\label{exemple}}

Il s'agit ici de donner un exemple de mesure $\nu=T \wedge S$, qui ne
charge pas d'ensemble algébrique et pour laquelle le plus petit
exposant $\chi_1(x)$ est nul pour presque tout point $x$.

Cette mesure va être construite à partir d'une application holomorphe
qui a un disque de Siegel. Plus précisément, soit $R(z)= \lambda z + z^2$ avec $\lambda=e^{2 i
  \pi \theta}$ où $\theta$ est irrationnel et diophantien (voir
\cite{CG} p.43 et p.55). $R$ a donc un disque de Siegel qui contient
$0$ (voir \cite{BM} p.33). En particulier $R$ est conjuguée à la
rotation $z \mapsto \lambda z$ au voisinage de $0$. C'est-à dire,
il existe une fonction holomorphe $\Phi$ localement inversible au
voisinage de $0$, pour laquelle $\Phi \circ R(z)= \lambda \Phi(z)$.

Si on considère un point $a$ de la courbe de
Jordan $|\Phi|= \epsilon$ (avec $\epsilon$ petit), alors $\frac{1}{m} \sum_{i=0}^{m-1}
\delta_{R^i(a)}$ converge vers la mesure de Lebesgue $\alpha$ de
$|\Phi|= \epsilon$. Cela provient du fait qu'une rotation d'angle
irrationnel sur le cercle unité est uniquement ergodique (voir
\cite{KH} p. 146). Maintenant, soit $f$ est l'endomorphisme holomorphe de
$\Pp^2(\Cc)$ défini par:
$$f([z:w:t])=[\lambda zt + z^2: \lambda wt+w^2:t^2].$$
Dans la suite, on se place dans la carte $t=1$.

L'image de la droite $z=a$ par $f^i$ est la droite $z=R^i(a)$
(parcourue $2^i$ fois). Autrement dit, si $S$ est une valeur
d'adhérence de $\frac{1}{m} \sum_{i=0}^{m-1} \frac{[f^i(L)]}{2^i}$,
alors $S= \int [N] d \alpha(N)$ où $\alpha$ est la mesure de Lebesgue précédente et $N$ est une droite qui passe par le point $[0:1:0]$ et
par un point de $| \Phi |= \epsilon$.

En particulier, comme $\alpha$ et $T \wedge [N]$ ne chargent aucun point, $T \wedge S$ ne
charge aucune courbe algébrique.

Enfin, comme $D$ est un disque de Siegel, il est assez facile de voir
que pour tout point $x$ du support de $\nu$ on a $\chi_1(x)=0$.

\newpage

\section{\bf{Appendice: la transformée de graphe}}

Dans cet appendice nous allons présenter une petite généralisation de
la transformée de graphe (voir \cite{KH}).

Le cadre de ce paragraphe est $\Cc^2$. Nous noterons $g$ l'application 
$$g(x,y)=(\lambda x + \alpha(x,y), \mu y + \beta(x,y))$$
avec $\alpha$, $\beta$ vérifiant $\alpha(0)= \beta(0)=0$, $| \alpha
|_{C^1} \mbox{ , } | \beta |_{C^1} \leq \delta$ dans une boule
$B=B(0,r) \subset \Cc^2$ et $0 < | \mu | < | \lambda |$. Soit maintenant $\{ (x, \phi(x)) \mbox{ , } x \in D \}$ un graphe
inclus dans $B$ au-dessus de $D$ qui vérifie $| \phi(x_1)-\phi(x_2) |
\leq \gamma |x_1-x_2 |$.

Le théorème de la transformée de graphe que nous allons énoncer donne
des conditions sur $\delta$, $\gamma$, $\lambda$ et $\mu$ pour que
l'image du graphe précédent par $g$ soit encore un graphe.

\begin{theoreme*}

Si $\delta (1 + \gamma) < | \lambda|$ alors l'image par $g$ du graphe
précédent est un graphe au-dessus de $\pi(g(\mbox{graphe de } \phi))$
(où $\pi$ est la projection sur l'axe des abscisses).

Par ailleurs, si $(x,\psi(x))$ désigne ce nouveau graphe, nous avons:
$$| \psi(x_1)- \psi(x_2) | \leq \frac{ | \mu | \gamma +
  \delta(1+\gamma)}{ | \lambda| - \delta(1 + \gamma)} |x_1-x_2 |.$$

\end{theoreme*}

\begin{proof}

La seule différence avec la démonstration qui se trouve dans \cite{KH}
p. $250$, c'est que nous ne sommes pas dans une situation
globale. Cependant, comme cela ne change rien à la preuve, nous ne
démontrerons que le premier point du théorème.

Soit $G: D \mapsto \Cc$ la fonction définie par $G(x)=\lambda x +
\alpha(x, \phi(x))$. C'est l'abscisse du point $g(x, \phi(x))$. Pour
voir que $g(\mbox{graphe de } \phi)$ est un graphe au-dessus de
$\pi(g(\mbox{graphe de } \phi))$, il suffit de voir que $G$ est une
bijection de $D$ sur $G(D)$. Autrement dit, si $x_0 \in G(D)$, il faut
démontrer que $G(x)=x_0$ (avec $x \in D$) a une unique solution.

Soit $F(x)= \lambda^{-1} x_0 -  \lambda^{-1} \alpha(x,
\phi(x))$. L'équation $G(x)=x_0$ équivaut à $F(x)=x$.
Par ailleurs, si $x_1 \mbox{, } x_2 \in D$, on a:
$$| F(x_1) - F(x_2) | = | \lambda ^{-1} \alpha(x_1, \phi(x_1)) -
\lambda^{-1} \alpha(x_2, \phi(x_2))| \leq | \lambda |^{-1} |
\alpha|_{C^1} d((x_1, \phi(x_1)),(x_2, \phi(x_2)))$$
car $(x_i, \phi(x_i)) \in B$ pour $i=1,2$.

On obtient donc
$$| F(x_1) - F(x_2) | \leq | \lambda |^{-1} \delta (1 + \gamma)
|x_1-x_2 |.$$
Maintenant si $\delta (1 + \gamma) < |\lambda|$ alors $| F(x_1) -
F(x_2) | \leq t |x_1-x_2 |$ avec $t < 1$ et l'équation $F(x)=x$ a donc
bien une unique solution dans $D$.

\end{proof}

\bigskip

Henry de Thélin\\
Université Paris-Sud (Paris 11)\\
Mathématique, Bât. 425\\
91405 Orsay\\
France

\end{document}